\documentclass[12pt,twoside,a4paper]{article}
\usepackage{graphicx}
\usepackage{amsfonts}
\usepackage{bbm}
\usepackage{mathrsfs}
\usepackage{fancyhdr}
 \usepackage{mathrsfs,amsmath,amssymb,amsthm}
\usepackage[dvips]{color}
 \usepackage{amssymb}
 \usepackage{amsbsy}
 \usepackage[dvips]{color}

\allowdisplaybreaks

 \theoremstyle{plain}
\theoremstyle{remark}  \newtheorem{remark}{\noindent\mbox{Remark}}
 \theoremstyle{plain}
 \theoremstyle{plain}\newtheorem{lemma}{\noindent\mbox{Lemma}}
\theoremstyle{plain} \newtheorem{theorem}{\noindent\mbox{Theorem}}
 \theoremstyle{plain}\newtheorem{proposition}{\noindent\mbox{Proposition}}
 \theoremstyle{plain}\newtheorem{corollary}{\noindent\mbox{Corollary}}
\theoremstyle{definition} 

 \def\proof{\noindent{\it Proof.~~}}
 \def\qed{\hfill$\Box$\medskip}
 \def\rto{\rightarrow\infty}
\def\z{\left}
\def\y{\right}
 \def\no{\nonumber}

\begin{document}
 \title{{Asymptotics of product of nonnegative 2-by-2 matrices  with applications to random walks with asymptotically zero drifts}}                

\author{Hua-Ming \uppercase{Wang}\footnote{Email:hming@ahnu.edu.cn; School of Mathematics and Statistics, Anhui Normal University, Wuhu, 241003, China }~ $\ \&\ $
Hongyan \uppercase{Sun}\footnote{Email:sun$\_$hy@cugb.edu.cn; Department of mathematics, China University of Geosciences, Beijing, 100083, China}
}
\date{}
\maketitle%

\vspace{-1cm}

\begin{center}
\begin{minipage}[c]{12cm}
\begin{center}\textbf{Abstract}\quad \end{center}
Let $A_kA_{k-1}\cdots A_1$ be product of some nonnegative 2-by-2 matrices. In general, its elements are hard to evaluate. Under some conditions, we show that  $\forall i,j\in\{1,2\},$ $(A_kA_{k-1}\cdots A_1)_{i,j}\sim c\varrho(A_k)\varrho(A_{k-1})\cdots \varrho(A_1)$ as $k\rightarrow\infty,$  where $\varrho(A_n)$ is the spectral radius of the matrix $A_n$ and $c\in(0,\infty)$ is some constant, so that the elements of $A_kA_{k-1}\cdots A_1$ can be estimated.
As applications, consider the maxima of certain excursions of (2,1) and (1,2) random walks with asymptotically zero drifts.
 We get some delicate limit theories which are quite different from the ones of simple random walks. Limit theories of both the tail and critical tail sequences of continued fractions play important roles in our studies.

\vspace{0.2cm}

\textbf{Keywords:}\  Product of nonnegative matrices, random walk,  spectral radius, tail of continued fraction.
\vspace{0.2cm}

\textbf{MSC 2010:}\ 15B48, 60G50, 60J10
\end{minipage}
\end{center}

\section{Introduction}
Product of nonnegative matrices are very useful in many aspects, for example, the multitype branching processes in random or varying environments(\cite{jlp,va}), random walk in random environment with bounded jumps (\cite{br}), Bernoulli convolutions and Gibbs properties of
linearly representable measures (\cite{ot}), etc.  As such, it has been extensively studied in literatures, see for example \cite{co,dl,fs,hj,235,jon} and  references therein. We also refer the reader to \cite{ha,se} for more general studies of  product of nonnegative matrices and their applications in studies of Markov chains.

In this paper, we study especially product $A_kA_{k-1}\cdots A_1,k\ge 1$ of 2-by-2 nonnegative matrices.
In applications, sometimes we need to evaluate the elements of $A_kA_{k-1}\cdots A_1,$ as $k\rto.$ For example, both the escape probabilities of nonhomogeneous (2,1) and (1,2) random walks from certain intervals can be written in terms of the elements of nonhomogeneous product of nonnegative 2-by-2 matrices, so that to study the limit behaviors of random walks, an inevitable step is to evaluate the elements of product of nonnegative 2-by-2 matrices.
Note that under some condition, it can be shown that $\frac{A_kA_{k-1}\cdots A_1}{\|A_kA_{k-1}\cdots A_1\|}$ or $\frac{\mathbf vA_kA_{k-1}\cdots A_1}{\|\mathbf vA_kA_{k-1}\cdots A_1\|}$ converge as $k\rto,$ where $\mathbf v$ is a row vector, see for example  \cite{fs,ot}. However, this fact provides no useful information for evaluating  elements of $A_kA_{k-1}\cdots A_1,$ since both $\|\mathbf vA_kA_{k-1}\cdots A_1\|$ and $\|A_kA_{k-1}\cdots A_1\|$ are hard to estimate in general.

As such, we turn to consider the spectral radius of $\varrho(A_n)$ of the matrix $A_n,n\ge1.$
Suppose that $A_k$ converges properly to a limit as $k\rto.$  We show that for $i,j\in\{1,2\},$ $\mathbf e_iA_kA_{k-1}\cdots A_1\mathbf e_j^t\sim c\varrho(A_k)\varrho(A_{k-1})\cdots \varrho(A_1),$ as $k\rto$ with $c\in(0,\infty)$ certain constant, $\mathbf e_1:=(1,0)$ and $\mathbf e_2:=(0,1).$ In this way, the elements of $A_kA_{k-1}\cdots A_1$ can be estimated.
As applications, we get some delicate limit theories of  the maxima of both (2,1) and (1,2) random walks, which are quite different from those of simple random walks.

\vspace{-.3cm}

\subsection{Product of nonnegative 2-by-2 matrices}
Let $a_k,b_k,d_k,$ $k\ge 1$ be certain positive numbers and for $k\ge1,$ set  \begin{align}
  A_k=\left(\begin{array}{cc}
                                                             a_k & b_k\\
                                                             d_k& 0
                                                           \end{array}
  \right).\label{ak}
\end{align}
Throughout the paper, for a matrix $A,$ we always denote by $\varrho(A)$ its spectral radius.
We now introduce some conditions on the numbers $a_k,b_k$ and $d_k, k\ge1.$

\noindent{\bf(B1)} For some $\sigma>0,$ $a_k, b_k,d_k \ge \sigma$ for all $k\ge1$ and
$$\sum_{k=2}^\infty|a_k-a_{k-1}|+|b_k-b_{k-1}|+|d_k-d_{k-1}|<\infty. $$

Notice that under condition (B1), there are some constants $0<c_1<c_2<\infty$ and $c_1\le a, b,d \le c_2$ such that \begin{align}
 &c_1\le a_k,b_k,d_k, \varrho(A_k) \le c_2 \text{ for all } k\ge 1 ,\label{ub}\\
 &\lim_{k\rto} a_k=a, \lim_{k\rto} b_k=b, \lim_{k\rto} d_k=d. \label{abc}
\end{align}
Thus under condition (B1), $A_k$ is convergent.

Suppose now condition (B1) holds. We introduce further the following conditions.

\noindent{\bf(B2)$_{\bf a}$} $\exists k_0>0,$ such that $\frac{a_k}{b_k}=\frac{a_{k+1}}{b_{k+1}},
\ \frac{d_k}{b_k}\ne\frac{d_{k+1}}{b_{k+1}},\ \forall k\ge k_0$ and
$$\lim_{k\rto}\frac{d_{k}/b_{k}-d_{k+1}/b_{k+1}}{d_{k+1}/b_{k+1}-d_{k+2}/b_{k+2}}$$
exists as an extending number.

 \noindent{\bf(B2)$_{\bf b}$} $\exists k_0>0,$ such that $\frac{a_k}{b_k}\ne\frac{a_{k+1}}{b_{k+1}},
\ \frac{d_k}{b_k}=\frac{d_{k+1}}{b_{k+1}},\ \forall k\ge k_0$ and
$$\lim_{k\rto}\frac{a_{k}/b_{k}-a_{k+1}/b_{k+1}}{a_{k+1}/b_{k+1}-a_{k+2}/b_{k+2}}$$ exists as an extending number.

\noindent{\bf(B2)$_{\bf c}$} $\exists k_0>0,$ such that $\frac{a_k}{b_k}\ne\frac{a_{k+1}}{b_{k+1}},
\ \frac{d_k}{b_k}\ne\frac{d_{k+1}}{b_{k+1}},\ \forall k\ge k_0$ and
   $$\tau:=\lim_{k\rto}\frac{d_{k}/b_k-d_{k+1}/b_{k+1}}{a_{k}/b_k-a_{k+1}/b_{k+1}}\ne \frac{-a+ \sqrt{a^2+4bd}}{2b}  $$ exists as an extending number. In addition, if $\tau$ is finite,  assume further $\lim_{k\rto}\frac{a_{k}/b_{k}-a_{k+1}/b_{k+1}}{a_{k+1}/b_{k+1}-a_{k+2}/b_{k+2}}$ exists as an extending number. Otherwise, if $\tau=\infty,$ assume further  $\lim_{k\rto}\frac{d_{k}/b_{k}-d_{k+1}/b_{k+1}}{d_{k+1}/b_{k+1}-d_{k+2}/b_{k+2}}$ exists as an extending number.
\begin{theorem}\label{amij} Suppose condition (B1) and one of (B2)$_a,$ (B2)$_b$ and  (B2)$_c$ hold. Then  $\forall i,j\in \{1,2\},$  for each $m\ge1,$
   there exists $0<c(m)<\infty$ such that
  \begin{align}\label{rac}
    \lim_{k\rto}\frac{\mathbf e_iA_k\cdots A_1\mathbf e_j^t}{\varrho(A_k)\cdots\varrho(A_1)}=c(m).
  \end{align}
  Moreover, we have $0<\varliminf_{m\rto}c(m)\le \varlimsup_{m\rto}c(m)<\infty.$
\end{theorem}
\begin{remark}\label{rem} (i)
  If $\exists k_0>0,$ such that $\frac{a_k}{b_k}=\frac{a_{k+1}}{b_{k+1}},
\ \frac{d_k}{b_k}=\frac{d_{k+1}}{b_{k+1}},\ \forall k\ge k_0,$ then $$A_{k}A_{k-1}\cdots A_{k_0}=a_ka_{k-1}\cdots a_{k_0}\z(\begin{array}{cc}
                                     1 & \lambda_1 \\
                                     \lambda_2 & 0
                                   \end{array}
\y)^{k-k_0+1}$$ for some $\lambda_1, \lambda_2>0,$ and (\ref{rac}) holds trivially.

(ii) The limits in conditions (B2)$_a,$ (B2)$_b$ and  (B2)$_c$ look awkward. Indeed, for example, if  both $\lim_{k\rto}\frac{d_k-d_{k+1}}{b_k-b_{k+1}}$ and $\lim_{k\rto}\frac{b_k-b_{k+1}}{b_{k+1}-b_{k+2}}$ exist as finite numbers,  then the limit in (B2)$_a$ exists. Roughly speaking,  it requires that  $a_k,$ $b_k$ and $d_k$ may fluctuate in different orders, but should fluctuate  in some common manner.
\end{remark}

 Let us  explain the main idea to prove Theorem \ref{amij}.  For $k\ge1,$ let $x_k:=\frac{\mathbf e_1A_k\cdots A_1\mathbf e_1^t}{\varrho(A_k)\cdots\varrho(A_1)}.$
 As the first step, under condition (B1), in Lemma \ref{lwb} below, we show that $x_k,k\ge1$ is uniformly bounded away from $0$ and infinity.  For this purpose, owing to  \cite{235},
 we show first in Lemma \ref{rfm} below that $\zeta\le \frac{\varrho(A_k\cdots A_1)}{\varrho(A_k)\cdots \varrho(A_1)}\le \gamma, \forall k\ge1$  for some $\zeta,\gamma\in(0,\infty).$ Then using the ergodiciy theorem of product of nonnegative matrices(\cite{se}) and the theory of limit periodic continued fractions,  we show in Lemma \ref{ral} that $\mathbf e_1A_k\cdots A_1\mathbf e_1^t \sim c\varrho(A_k\cdots A_1)$ for some $c\in (0,\infty)$ as $k\rto.$ Therefore for some $c_3,c_4\in(0,\infty),$   $\forall k\ge1,$ we have $c_3\le x_k\le c_4.$

  As the second step,  we show that $x_k-x_{k-1}$ converges to $0$ either exponentially fast or in an alternating manner. To this end, we first develop some fluctuation theory of the {\it critical tail sequence} of a limit periodic continued fraction, where one of conditions (B2)$_a,$ (B2)$_b$ and  (B2)$_c$ is required, see Lemmas \ref{de} and \ref{eq}  below. Then, the fluctuation of $x_k-x_{k-1}$ can be studied.  As a result,   $x_k$ converges to some number $c\in(0,\infty)$ as $k\rto.$

\vspace{-.3cm}

\subsection{Maxima of (2,1) and (1,2) random walks}
Next, as applications of Theorem \ref{amij}, we consider the maxima of certain excursions of (2,1) and (1,2) random walks with asymptotically zero drift. Firstly, let us introduce precisely the model.
Suppose that \begin{align*}
  q_k,p_k,k\ge 2\text{ are numbers such that }\forall k\ge2, q_k,p_k>0,
  \ q_k+p_k=1.
\end{align*} Let $Y=\{Y_k\}_{k\ge 0}$ be a Markov chain on $\mathbb Z_+:=\{0,1,2,...\}$ starting from some $y_0\in\mathbb Z_+$ and with transition probabilities
\begin{align}
&P(Y_{k+1}=1|Y_k=0)=P(Y_{k+1}=2|Y_k=1) =1\nonumber,\\
&P(Y_{k+1}=n+1|Y_k=n)= q_n,\nonumber\\
&P(Y_{k+1}=n-2|Y_k=n)= p_n, n\ge 2,k\ge0.\nonumber
\end{align}
Introduce also another Markov chain $Y'=\{Y'_k\}_{k\ge0}$ on $\mathbb Z_+,$ starting from some $y_0'\in \mathbb Z_+$ and with transition probabilities
\begin{align*}
 &P(Y'_{k+1}=0|Y'_k=1)=P(Y'_{k+1}=2|Y'_k=0)=1,\\
  &P(Y'_{k+1}=n-1|Y'_k=n)=q_n,\\
  &P(Y'_{k+1}=n+2|Y'_k=n)= p_n,n\ge2,k\ge 0.
\end{align*}

Unless otherwise stated, we always assume that both $Y$ and $Y'$ start from $y_0=y_0'=2.$ We call the chain $Y$ a (2,1) random walk and $Y'$ a (1,2) random walk.
In literatures,  $Y'$ is called the adjoint chain of $Y$ and vice versa.

Next, we define the so-called maxima we concern for the chains $Y$ and $Y'.$  For $X \in \{Y,Y'\},$ we denote
\begin{equation*}\begin{split}
  &D(X):=\inf\{k\ge 1: X_k<X_0\},\\
  & M(X):=\sup\{X_k:0\le k\le D(X)\},\end{split}
\end{equation*}
where and throughout, we assume $\inf \phi=\infty.$ Clearly, $D(X)$ is
the time that the chain $X$ hits some point below $X_0$ for the first time and $M(X)$ is the maximum of the excursion $\{X_0,X_1,...,X_D\}.$
When no danger of making confusion, we write $D(X)$ and $M(X)$ as $D$ and $M$ respectively for simplicity.

To begin with, let us  consider  (2,1) random walk $Y.$
We aim to study the distribution of $M(Y)$ and  characterize its asymptotics as well. To this end,
for $k\ge 2,$ introduce matrix
\begin{align}
N_k:=\z(\begin{array}{cc}
          \theta_k & \theta_k \\
          1 & 0
        \end{array}
\y)\text{ with }\theta_k:=\frac{p_k}{q_k},\label{tn}
\end{align}
which we will work with.
Proposition \ref{dm} below gives the distribution of $M(Y)$ on the event $\{D(Y)<\infty\}$ in terms of $N_k,k\ge2.$
\begin{proposition}\label{dm} Consider (2,1) random walk $Y.$ For $n\ge 2,$ we have
  \begin{align}
    P(&M=n,D<\infty)
    =\frac{1}{1+\sum_{s=2}^{n-1}\mathbf e_1 N_s\cdots N_{2}\mathbf e_1^t}\frac{\mathbf e_1 N_n\cdots N_{2}\mathbf e_1^t}{1+\sum_{s=2}^{n}\mathbf e_1 N_s\cdots N_{2}\mathbf e_1^t}.\label{pmn}
  \end{align}
\end{proposition}
Firstly, let us see what happens to  simple random walk. For $k\ge 2,$ set
\begin{align*} q_k\equiv q\in (0,1), p_k\equiv p=1-q\text{ and  } N_k\equiv N:= \z(\begin{array}{cc}
   p/q & p/q \\
1 & 0
 \end{array}
\y).
\end{align*}
In this case, the chain $Y$ is positive recurrent, null recurrent or transient according as $q<2/3, =2/3\text{ or }>2/3$(or equivalently, $\varrho(N)>1,=1 \text{ or } <1$),  respectively. Some direct computation from (\ref{pmn}) yields that
$$P(M(Y)=n,D(Y)<\infty)\sim \z\{ \begin{array}{cc}
                             c\varrho(N)^{-n},  & \text{if }\varrho(N)>1, \\
                             c/n^2, & \text{if }\varrho(N)=1, \\
                             c\varrho(N)^n, & \text{if }\varrho(N)<1,
                           \end{array}
\y. \text{ as }n\rto.$$
Here and in the rest of the paper, unless otherwise specified, $0<c<\infty$ is some constant, whose value may change from line to line.

We conclude that $P(M(Y)=n,D(Y)<\infty)$ decays either {\it exponentially} if the walk is transient or positive recurrent, or {\it polynomially} with speed $c/n^2$ if the walk is null recurrent.

We next study  (2,1) and (1,2) random walks with asymptotically zero drifts, for which, $P(M=n,D<\infty)$ decays with  various speeds quite different from those of simple random walk.
Adding some perturbations on the transition probabilities of a null recurrent simple random walk, we get a near-recurrent random walk, known also as Lamperti random walk which dates back to Harris \cite{har} and Lamperti \cite{lam} and has been extensively studied in literatures, refer for example, to \cite{cfr,cfrb,hmw13,jlp,lo,w19} etc. To introduce Lamperti random walk,
we take the perturbations from \cite{cfr}.
For $K= 1,2,...$ and $B\in \mathbb R,$ set
\begin{align*}
    &\Lambda(1,i,B)=\frac{B}{i},\\
   &\Lambda(2,i,B)=\frac{1}{i}+\frac{B}{i\log i}, \cdots,\\
    &\Lambda(K,i,B)=\frac{1}{i}+\frac{1}{i\log i}+...+\frac{1}{i\log i\cdots\log_{K-2}i}+\frac{B}{i\log i\cdots\log_{K-1}i},
  \end{align*}
  where $\log_0 i=i$ and for $k\ge1,$ $\log_{k}i=\log \log_{k-1} i.$

 For $K$ and $B$ fixed, set $$i_0:=\min\z\{i: \log_{K-1}i>0, {|\Lambda(K,i,B)|}< 1 \y\},$$
 and let
\begin{align}r_i:=\left\{\begin{array}{ll}
   \frac{\Lambda(K,i,B)}{3}, & i\ge  i_0, \\
 r_{i_0}, & i< i_0,
\end{array}\right.\label{ri}\end{align}
which serves as perturbations added on a null recurrent simple random walk.
\begin{theorem}\label{mxt} Consider (2,1) random walk $Y.$  Fix $K=1,2,3,...$ and $B\in \mathbb R.$
 {\rm (i)} If
    $q_i=\frac{2}{3}+r_i,i\ge 2,$
then, as $n\rto,$
\begin{align}\label{ap}
  P(M=n,D<\infty)\sim \left\{ \begin{array}{ll}
                                \frac{c }{n\log n\cdots \log_{K-2}n\log_{K-1} n(\log_K n)^2}, &  \text{if } B=1,  \\
\frac{c}{n\log n\cdots \log_{K-2}n(\log_{K-1} n)^B}, &  \text{if } B>1,  \\
\frac{c}{n\log n\cdots \log_{K-2}n(\log_{K-1} n)^{2-B}}, &  \text{if } B<1.
                              \end{array}
  \right.
\end{align}
{\rm(ii)} If
$  q_i=\frac{2}{3}-r_i,i\ge 2,$
then, as $n\rto,$
\begin{align}\label{an}
  P(M=n,D<\infty)\sim
  \left\{\begin{array}{ll}
\frac{c}{n^{B+2}}, &\text{if } K=1,B>-1,\\
\frac{c}{n(\log n)^2}, & \text{if }K=1, B=-1,\\
  cn^B, & \text{if }K=1,B<-1,\\
 \frac{c}{n^3\log n...\log_{K-2}n (\log_{K-1} n)^B}, &\text{if } K>1.\\
\end{array}\right.
\end{align}
\end{theorem}
\begin{remark} {\bf (a)} If $K=1,$
since  $\Lambda(1,i,B)=\frac{B}{i},$
 replacing $B$ by $-B$ in (\ref{ap}), we get \eqref{an} and vice versa.

 {\bf(b)} Fixing $K=1$ and letting $q_i=2/3+r_i,i\ge 2,$ we shed some light on the null recurrent case to illustrate the difference between Lamperti random walk and simple random walk. If $-1\le B\le 1,$ then by Corollary \ref{crpt} below, the chain $Y$ is null recurrent
  and thus $P(D<\infty)=1.$ By \eqref{ap}, we have
   \begin{align*}
  P(M=n)\sim \left\{ \begin{array}{ll}
                                \frac{c }{ n(\log n)^2}, &  \text{if } B=1,  \\
\frac{c}{n^{2-B}}, &  \text{if } -1\le B<1,
                              \end{array}
  \right. \text{ as }n\rto.
\end{align*}
  So, even for null recurrent case, the decay speeds are quite  sensitive to $B.$ But, for null recurrent simple random walk, $P(M=n)$ always decays polynomially with speed $cn^{-2}.$

  {\bf(c)} The difficulty for (2,1) random walk arises from  the fact that the escape probabilities are functions of the product of matrices, $\mathbf e_1N_k\cdots N_2\mathbf e_1^t,$ which cannot be estimated directly.
  However, by Theorem \ref{amij} we have
  \begin{align}\mathbf e_1N_k\cdots N_2\mathbf e_1^t\sim c\varrho(N_k)\cdots \varrho(N_2)\text{ as }k\rto.\label{srpm}\end{align}
  Therefore
  it suffices to work with
   $\varrho(N_k)\cdots \varrho(N_2),$ which  is able to be estimated.

  {\bf(d)}  For the nearest neighbor random walk, similar result holds, see \cite{wh}.
  In this setting, if we denote by $p'_k$ the probability that the walk jumps from $k$ to $k+1$ in the next step whenever it is currently located at $k,$  then with $\rho_k:=\frac{1-p'_k}{p'_k},$ the escape probabilities are written in terms of $\rho_1\cdots\rho_n,$ which can be directly estimated. Therefore for this setting, things are much easier.
  \end{remark}
Next, we consider (1,2) random walk $Y'.$ To derive similar result, besides the asymptotics of product of nonnegative matrices, we need to develop further some other techniques related to the limit periodic continued fractions  and the hitting probabilities of the walk.
\begin{theorem}\label{myp} Consider (1,2) random walk $Y'.$
 Fix $K=1,2,3,...$ and $B\in \mathbb R.$
 {\rm (i)} If
    $p_i=\frac{1}{3}+r_i,i\ge 2,$
then, \eqref{ap} holds as $n\rto.$
{\rm(ii)} If
$  p_i=\frac{1}{3}-r_i,i\ge 2,$
then, \eqref{an} holds as $n\rto.$
\end{theorem}
We now explain the idea and difficulty to prove Theorem \ref{myp}.
For  $1\le m\le k\le n,$ $j\in \{n,n+1\},$  let
\begin{align*}
&\mathcal P_k^j(m,n,+)=P(Y'\text{ hits }[n,\infty] \text{ at }j\text{ before }[0,m] |Y'_0=k),\\
&\mathcal P_k(m,n,+)=P(Y'\text{ hits }[n,\infty] \text{ before }[0,m] |Y'_0=k).
\end{align*}  Then by Markov property, we can get \begin{equation}P(M=n,D<\infty)=\mathcal P_2^n(1,n,+)(1-\mathcal P_n(1,n+1,+)).\label{yp}\end{equation}
Similar to Lemma \ref{esy} below, by some computation, it can be shown that
\begin{equation}\label{ps}\begin{split}
  &1-\mathcal P_n(1,n+1,+)=\frac{1}{1+\sum_{s=2}^{n}\mathbf e_1N_s\cdots N_{n}\mathbf e_1^t},\\
  &\mathcal P_2^n(1,n,+)=\mathbf e_1N_2\cdots N_{n-1}\Big(\frac{1+\sum_{s=2}^{n-1}\mathbf e_1N_s\cdots N_{n-1}\mathbf e_2^t}{1+\sum_{s=2}^{n-1}\mathbf e_1N_s\cdots N_{n-1}\mathbf e_1^t}\mathbf e_1^t-\mathbf e_2^t\Big),\end{split}
\end{equation}
which are hard to estimate even though we have \eqref{srpm} in hands, since every summand there depends on  $n$.
The good news is that, in \cite{leta},  the escape probability $\mathcal P_2(1,n,+)$
(NOT $\mathcal P_2^n(1,n,+)$ ) can be written in terms of the tail $\xi_n:=\frac{\theta_{n}^{-1}}{1}\begin{array}{c}
                                \\
                               +
                             \end{array}\frac{\theta_{n+1}^{-1}}{1}\begin{array}{c}
                                \\
                               +\cdots
                             \end{array}$ of a continued fraction, whose estimation has been given in \cite{w19}. However, to study the asymptotics of $P(M=n,D<\infty),$ we need to know exactly what $\mathcal P_2^n(1,n,+)$ is, or at least, to find out its asymptotics. By constructing a new Markov chain related to $Y'$ and analyzing the hitting probabilities of the new chain, Lemma \ref{rnn} below shows that $\frac{\mathcal P_2^n(1,n,+)}{\mathcal P_2^{n+1}(1,n,+)}\rightarrow 2,$  so that the term $\mathcal P_2^n(1,n,+)$ in \eqref{yp} can be estimated.

The formula for $1-\mathcal P_n(1,n+1,+)$ in \eqref{ps} looks somewhat simple, but some special techniques are required to give its estimation. In Lemma \ref{mc} below, with the help of \eqref{srpm}, by a delicate analysis of the product $N_s\cdots N_n$ and the continued fraction, we show that
$1-\mathcal P_n(1,n+1,+)\sim c\frac{\xi_2\cdots\xi_n}{\sum_{s=2}^{n+1}\xi_2\cdots \xi_{s-1}} $  as $n\rto.$ So Theorem \ref{myp} can finally be proved.


The rest of the paper is organized as follows. Section \ref{spm} is devoted to proving Theorem \ref{amij}. In Section \ref{mrw}, we study the (2,1) random walk $Y$ to give the proof of Theorem \ref{mxt}. The maximum for (1,2) random walk $Y'$ are studied in Section \ref{mt}.
Finally, the criteria for transience, null recurrence and positive recurrence of both the chains $Y$ and $Y'$ are presented in an appendix section.
\section{Asymptotics of  product of nonnegative 2-by-2 matrices}\label{spm}
The main task of this section is to prove Theorem \ref{amij}.
To begin with, we introduce some notations of continued fraction which will be used time and time again.

Let $\beta_k,\alpha_k,k\ge 1$ be certain  numbers. We denote by
\begin{equation*}
\mathrm K_{n=1}^\infty(\beta_n|\alpha_n)\equiv\frac{\beta_1}{\alpha_1}\begin{array}{c}
                                \\
                               +
                             \end{array}\frac{\beta_2}{\alpha_2 }\begin{array}{c}
                                \\
                               +
                             \end{array}\frac{\beta_3}{ \alpha_3}\begin{array}{c}
                                \\
                               +\cdots
                             \end{array}:=\dfrac{\beta_1}{\alpha_1+\dfrac{\beta_2}{{\alpha_2+\dots}}}
\end{equation*}
 a continued fraction, and by
 \begin{align}
  &f^{(n)}:=\frac{\beta_{n+1}}{\alpha_{n+1}}\begin{array}{c}
                                \\
                               +
                             \end{array}\frac{\beta_{n+2}}{\alpha_{n+2 }}\begin{array}{c}
                                \\
                               +
                             \end{array}\frac{\beta_{n+3}}{\alpha_{n+3}}\begin{array}{c}
                                \\
                               +\cdots
                             \end{array},\ n\ge 0,\no\\
&h_k:=\frac{\beta_k}{\alpha_{k-1}}\begin{array}{c}
                                \\
                               +
                             \end{array}\frac{\beta_{k-1}}{\alpha_{k-2}}\begin{array}{c}
                                \\
                               +\cdots  +
                             \end{array}\frac{\beta_{2}}{\alpha_1}\begin{array}{c}
                             \end{array},\ k\ge 2 \no
\end{align}
 its $n$-th tail and its critical tail sequence respectively.

\subsection{Lower and upper bounds of $\frac{\mathbf{e}_1 A_{k} \cdots A_{1}  \mathbf e_1^t}{\varrho(A_{k}) \cdots \varrho(A_{1})}$ }

For $k\ge 1,$ let $A_k$ be the one in (\ref{ak}). We have
  \begin{align}
    \varrho(A_k)=\frac{a_k+\sqrt{a_k^2+4b_kd_k}}{2}.\label{rhk}
  \end{align}
In this subsection, we show that $\frac{\mathbf{e}_1 A_{k} \cdots A_{1}  \mathbf e_1^t}{\varrho(A_{k}) \cdots \varrho(A_{1})},k\ge1$ is uniformly bounded away from $0$ and infinity.
\begin{lemma}\label{lwb}
  Suppose that  condition (B1) holds. Then there exist constants $0<c_3<c_4<\infty$ such that  for  $k\ge m\ge1,$ $c_3<\frac{\mathbf{e}_1 A_{k} \cdots A_{m}  \mathbf e_1^t}{\varrho(A_{k}) \cdots \varrho(A_{m})}<c_4.$
\end{lemma}
\proof The lemma is a direct consequence of Lemmas \ref{rfm} and  \ref{ral} below.\qed
   \begin{lemma}\label{rfm}
    Suppose that  condition (B1) holds.  Then there exist constants $0<\zeta<\gamma<\infty$ such that for $k\ge m\ge1,$
    $
 \zeta \leq \frac {\varrho(A_{k} \cdots A_{m})}{\varrho(A_{k}) \cdots \varrho(A_{m})}\leq \gamma.
$

  \end{lemma}

\proof For vectors $\mathbf v= \left(
        \begin{array}{c}
        v_1\\
        v_2  \\
        \end{array}
       \right)$ and $\mathbf u=\left(
        \begin{array}{c}
        u_1\\
        u_2  \\
        \end{array}
       \right),$ set $$\frac{\mathbf v}{\mathbf u}:= \left(
        \begin{array}{c}
        v_1/u_1\\
        v_2/u_2  \\
        \end{array}
       \right),\ \left(\frac{\mathbf v}{\mathbf u}\right)_{\textrm{min}}:=\min\z\{\frac{v_1}{u_1},\frac{v_2}{u_2}\y\}, \ \left(\frac{\mathbf v}{\mathbf u}\right)_{\textrm{max}}:=\max\z\{\frac{v_1}{u_1},\frac{v_2}{u_2}\y\}.$$
  Let $\mathbf v_n$ be a right eigenvector of $A_n$ corresponding to the largest eigenvalue $\varrho(A_n)$. Then we can choose $\mathbf v_n$ to be
$$
\mathbf{v}_n=\left(
        \begin{array}{c}
        \varrho(A_n)\\
        d_n  \\
        \end{array}
       \right).
$$
For $k\ge m\ge 1,$ write \begin{align}
  &\gamma_{k,m}:=\Big(\frac{\mathbf v_k}{\mathbf v_{k-1}}\Big)_{\max} \cdots \Big(\frac{\mathbf v_{m+1}}{\mathbf v_{m}}\Big)_{\max}
\Big(\frac{\mathbf v_{m}}{\mathbf v_{k}}\Big)_{\max},\label{ga}\\
&\zeta_{k,m}:=\Big(\frac{\mathbf v_k}{\mathbf v_{k-1}}\Big)_{\min} \cdots \Big(\frac{\mathbf v_{m+1}}{\mathbf v_{m}}\Big)_{\min}
\Big(\frac{\mathbf v_{m}}{\mathbf v_{k}}\Big)_{\min}.\label{ze}
\end{align}
Applying \cite[Theorem 1, page 228]{235}, for $k\ge1,$ we have
\begin{align*}
\zeta_{k,m}\leq
\frac {\varrho(A_{k} \cdots A_{m})}{\varrho(A_{k}) \cdots \varrho(A_{m})}\le \gamma_{k,m}.
\end{align*}
It remains to show that both $\zeta_{k,m}^{-1}$ and $\gamma_{k,m},k\ge m\ge1$ are uniformly bounded away from $\infty.$  To this end, set $\epsilon_n=(\mathbf v_n/\mathbf v_{n-1})_{\max}-1, n\ge2.$ Then by (\ref{ub}),
\begin{align*}
  |\epsilon_n|&\le\max\left\{\left|\frac{\varrho(A_n)}{\varrho(A_{n-1})}-1\right|, \left|\frac{d_n}{d_{n-1}}-1\right|\right\}
  \le \left|\frac{\varrho(A_n)}{\varrho(A_{n-1})}-1\right|+\left|\frac{d_n}{d_{n-1}}-1\right|\\
  &\le c(|\varrho(A_{n})-\varrho(A_{n-1})|+|d_{n}-d_{n-1}|)\\
  &\le c(|a_n-a_{n-1}|+|b_n-b_{n-1}|+|d_n-d_{n-1}|), n\ge 2.
\end{align*}
Therefore, by condition (B1), $\sum_{n=2}^{\infty}|\epsilon_n|<\infty, $ implying that $\sum_{n=2}^\infty \log(1+|\epsilon_n|)<\infty.$ As a consequence, taking \eqref{ub} and \eqref{abc} into account, we have $$\gamma_{k,m}\le \max\Big\{\frac{\varrho(A_m)}{\varrho(A_k)},
\frac{d_m}{d_k}\Big\}\prod_{n=2}^\infty(1+|\epsilon_n|)<\gamma$$ for some number $\gamma<\infty$ independent of $k$ and $m.$

Since $\zeta_{k,m}^{-1}=\Big(\frac{\mathbf v_{k-1}}{\mathbf v_k}\Big)_{\max} \cdots \Big(\frac{\mathbf v_{m}}{\mathbf v_{m+1}}\Big)_{\max}
\Big(\frac{\mathbf v_{k}}{\mathbf v_{m}}\Big)_{\max},$ a similar argument also yields that
$\zeta_{k,m}^{-1},k\ge m\ge1$ is uniformly bounded from above. Consequently, Lemma \ref{rfm} is proved.\qed
\begin{remark}\label{brb}
  If  condition (B1) is replaced by ``(B1)$'$: for some numbers $a,b,d>0,$  $a_k\rightarrow a,$ $b_k\rightarrow b,$  $d_k\rightarrow d$  as $k\rto,$ and in addition, $\frac{a_k}{d_k}, \frac{b_k}{d_k}$ are increasing (or decreasing) in $k$ simultaneously", then we may choose a right eigenvector corresponding to  $\varrho(A_n)$ as
$
\mathbf{v}_n=\left(
        \varrho(A_n)/d_n,
        1  \right)^t.
$ Let $\gamma_k$ and $\zeta_k$ be those defined in \eqref{ga} and \eqref{ze}. Under condition (B1)$'$, $\varrho(A_n)/d_n$ is monotone in $n.$ If it is increasing in $n,$ then
$\gamma_k=\frac{\varrho(A_k)d_1}{\varrho(A_1)d_k}$ and $\zeta_k=\frac{\varrho(A_1)d_k}{\varrho(A_k)d_1},k\ge 1.$
Since $\lim_{k\rto}\gamma_k=\lim_{k\rto}\zeta_k^{-1}=c$ for some $0<c<\infty,$ then both $\gamma_k,k\ge1$ and $\zeta_k,k\ge1$ are uniformly bounded away from $0$ and infinity.
Otherwise, if $\varrho(A_n)/d_n$ is decreasing  in $n,$  things can be done by a similar approach.
\end{remark}

  \begin{lemma} \label{ral} Suppose that  condition (B1) holds. For $n\ge 1,$ set $\tilde a_n=\frac{a_n}{b_n}, \tilde d_n=\frac{d_n}{b_n}.$  Let
    $$
g_m=\frac{ 1}{\tilde a_m}\begin{array}{c}
                                \\
                               +
                             \end{array}\frac{\tilde d_m}{\tilde a_{m+1}}\begin{array}{c}
                                \\
                               +
                             \end{array}\frac{\tilde d_{m+2}}{\tilde a_{m+3}}\begin{array}{c}
                                \\
                               +\cdots
                             \end{array},m\ge1\text{ and }f=\frac{\sqrt{a^2+4bd}-a}{2b},
$$ where $a,b,d$ are those numbers in \eqref{abc}.
Then we have $\lim_{m\rto}g_m=\frac{2b}{a+\sqrt{a^2+4bd}}$ and for each $m\ge1,$ $0<g_m<\infty$ and
\begin{align}
\lim_{k\rto}\frac{\varrho(A_{k} \cdots A_{m})}{\mathbf{e}_1 A_{k} \cdots A_{m}  \mathbf e_1^t}=1+fg_m. \label{kmp}
 \end{align}
  \end{lemma}

\proof
  For $k\ge 1,$ write
  $$A_k\cdots A_m=\z( \begin{array}{cc}
                        M_{k,m}(11) & M_{k,m}(12) \\
                        M_{k,m}(21)& M_{k,m}(22)
                      \end{array}
  \y).$$
    We claim that
  \begin{align}
  f&=\lim_{k\rto}\frac{M_{k,m}(21)}{M_{k,m}(11)}=\lim_{k\rto}\frac{M_{k,m}(22)}{M_{k,m}(12)},\label{fm}\\
  g_m&= \lim_{k\rto}\frac{M_{k,m}(12)}{M_{k,m}(11)}=\lim_{k\rto}\frac{M_{k,m}(22)}{M_{k,m}(21)}.\label{gm}
\end{align}
Indeed, an application of the ergodicity theorem of the product of nonnegative matrices(see \cite[Theorem 3.3]{se}) yields the existence of the limits  and the second equality in both (\ref{fm}) and (\ref{gm}).
 To compute $ f$ and $ g_m,$
noticing that
\begin{align*}
A_{k} \cdots A_{m}&=\left(
      \begin{array}{cc}
      a_{k} & b_{k} \\
      d_k & 0 \\
      \end{array}
    \right) \cdots \left(
      \begin{array}{cc}
      a_{m} & b_{m} \\
      d_m & 0 \\
      \end{array}
    \right)\\&=b_k \cdots b_{m} \left (
    \begin{array}{cc}
      \tilde a_k & 1 \\
     \tilde d_k & 0 \\
      \end{array}
    \right) \cdots \left(
      \begin{array}{cc}
      \tilde a_m & 1 \\
      \tilde d_m & 0 \\
      \end{array}
    \right),
    \end{align*}
thus  we have
\begin{align}
&\frac{M_{k,m}(21)}{M_{k,m}(11)}=\frac{\tilde d_{k}}{\tilde a_k}\underset{+}\quad \frac{\tilde d_{k-1}}{\tilde a_{k-1}}\underset{+}\quad \underset{\cdots}\quad \frac{\tilde d_{m+1}}{\tilde a_{m+1}}\underset{+}\quad \frac{\tilde d_{m}}{\tilde a_m}\label{cfa}=:f_{k,m},\\
&\frac{M_{k,m}(12)}{M_{k,m}(11)}=\frac{1}{\tilde a_m}\underset{+}\quad \frac{\tilde d_{m}}{\tilde a_{m+1}}\underset{+}\quad \frac{\tilde d_{m+1}}{\tilde a_{m+2}}\underset{+}\quad \underset{\cdots}\quad \frac{\tilde d_{k-2}}{\tilde a_{k-1}}\underset{+}\quad \frac{\tilde d_{k-1}}{\tilde a_k}=:g_{k,m}\label{cfb}
\end{align}
where \eqref{cfa} follows by forward induction and \eqref{cfb}  by backward induction.
It follows from (\ref{cfa}) that  \begin{align}\label{rk}f_{k,m}=\frac{\tilde d_k}{\tilde a_k+f_{k-1,m}},k\ge2.\end{align}
Since $\lim_{k\rto}\tilde a_k=a/b, \lim_{k\rto} \tilde d_k=d/b,$ letting $k\rto$ in \eqref{rk}, we get
$f=\frac{d/b}{a/b+ f},$ whose positive solution is
\begin{align}\label{ff} f=\frac{\sqrt{a^2+4bd}-a}{2b}.\end{align}
By the theory of  convergence of limit periodic continued fractions(see \cite[Theorem 4.13, page 188]{lw}), from (\ref{cfb}), we have
 \begin{align*}
 g_m&=\lim_{k\rto}g_{k,m}=\frac{ 1}{\tilde a_m}\begin{array}{c}
                                \\
                               +
                             \end{array}\frac{\tilde d_m}{\tilde a_{m+1}}\begin{array}{c}
                                \\
                               +
                             \end{array}\frac{\tilde d_{m+1}}{\tilde a_{m+2}}\begin{array}{c}
                                \\
                               +\cdots
                             \end{array} \in(0,\infty),\\
     g&:=\lim_{n\rto}g_m= \frac{2b}{a+\sqrt{a^2+4bd}}.
\end{align*}

Finally, we turn to prove \eqref{kmp}.
By some easy computation, we obtain
\begin{align}
  &\varrho(A_k\cdots A_m)=\frac{M_{k,m}(11)+M_{k,m}(22)}{2}\no\\
  &+\frac{\sqrt{(M_{k,m}(11)+M_{k,m}(22))^2+4(M_{k,m}(12)M_{k,m}(21)-M_{k,m}(11)M_{k,m}(22))}}{2}.\no
\end{align}
 Consequently, it follows from (\ref{fm}) and (\ref{gm}) that\begin{align*}
\lim_{k\rto}\frac{\varrho(A_{k} \cdots A_{m})}{\mathbf{e}_1 A_{k} \cdots A_{m}  \mathbf e_1^t}=\lim_{k \rto}\frac{\varrho(A_{k} \cdots A_{m})}{M_{k,m}(11)}=1+fg_m.
 \end{align*}
 The lemma is proved. \qed
\subsection{Critical tail sequence of a continued fraction}
In this subsection, we study the critical tail sequence of  a limit periodic continued fraction, which is helpful to prove Theorem \ref{amij}.

 \begin{lemma}\label{de} Suppose that $\alpha_k,\beta_k, \omega_k>0, k\ge1$ are numbers such that  $\alpha_k\rightarrow \alpha,$ $\beta_k\rightarrow \beta,$ and $\omega_k\rightarrow \frac{\sqrt{\alpha^2+4\beta}-\alpha}{2}=:\omega,$ as $k\rto,$ where $0<\alpha,\beta<\infty$ are certain constants. For $k\ge1,$ let
 \begin{align}f_k:=\frac{\beta_{k}}{\alpha_k}\underset{+}\quad \frac{\beta_{k-1}}{\alpha_{k-1}}\underset{+}\quad \underset{\cdots}\quad \frac{\beta_{2}}{\alpha_2}\underset{+}\quad \frac{\beta_{1}}{\alpha_1}.\label{ctf}
 \end{align}
 Set \begin{align}\label{dlep}\varepsilon_k=f_k-\omega_k, k\ge 1 \text{ and }\ \delta_{k}=\beta_k-\omega_{k}(\alpha_k+\omega_{k-1}),k\ge2.\end{align}
   Let $q$ be a fixed number. We have
  \begin{align*}
  &\text{if }\lim_{k\rto}\frac{\varepsilon_k}{\varepsilon_{k+1}}=q,\text{ then } |q|\ge1 \text{ and }\lim_{k\rto}\frac{\delta_k}{\delta_{k+1}}=q;\\
    &\text{if }\lim_{k\rto}\frac{\delta_k}{\delta_{k+1}}=q,\text{ then } |q|\ge1 \text{ and } \lim_{k\rto}\frac{\varepsilon_k}{\varepsilon_{k+1}}=q\text{ or }-\frac{1+\omega}{\omega}.
\end{align*}
 \end{lemma}
 \begin{remark}
   Similar results hold for the tail $f^{(n)}$ of  a limit periodic continued fraction, see \cite{jw90} and \cite{lor}. But for the critical tail sequence, up to our knowledge, there is no such observation yet.
 \end{remark}
\proof To prove the lemma, we follow \cite[Theorem 6.1, page 91]{lor}. It follows from (\ref{ctf}) that for  $k\ge 2,$ $f_k(\alpha_k+f_{k-1})=\beta_k$ and thus $$\delta_k=\varepsilon_{k-1}\omega_k+\varepsilon_{k}(\alpha_k+\omega_{k-1}+\varepsilon_{k-1}).$$
Therefore we get
\begin{align}\label{frde}
  \frac{\delta_k}{\delta_{k+1}}=\frac{\frac{\varepsilon_{k-1}}{\varepsilon_{k}}\omega_k+ \alpha_k+\omega_{k-1}+\varepsilon_{k-1}}{\frac{\varepsilon_k}{\varepsilon_{k+1}}\omega_{k+1}+\alpha_{k+1}
  +\omega_{k}+\varepsilon_{k}}\frac{\varepsilon_k}{\varepsilon_{k+1}}.
\end{align}
We first suppose that $\lim_{k\rto}\frac{\varepsilon_k}{\varepsilon_{k+1}}=q.$
By (\ref{fm}), (\ref{cfa}) and (\ref{ff}), we have $\lim_{k\rto}f_k=\frac{\sqrt{\alpha^2+4\beta}-\alpha}{2}=\omega,$ implying that $\lim_{k\rto}\varepsilon_k=0.$ Therefore $|q|\ge 1.$ Letting $k\rto$ in (\ref{frde}), we get $\lim_{k\rto} \frac{\delta_k}{\delta_{k+1}}=q.$

Next, assume  $\lim_{k\rto} \frac{\delta_k}{\delta_{k+1}}=q.$  Since $\lim_{k\rto}\delta_k=0,$ $|q|\ge 1.$ We can deduce from (\ref{frde}) that
\begin{align*}
  \frac{\varepsilon_k}{\varepsilon_{k+1}}\omega_{k+1}=
  \frac{\frac{\delta_k}{\delta_{k+1}}(\alpha_{k+1}
  +\omega_{k}+\varepsilon_{k})\omega_{k+1}}{\alpha_k+\omega_{k-1}+\varepsilon_{k-1}
  -\frac{\delta_k}{\delta_{k+1}}\omega_{k+1}+\frac{\varepsilon_{k-1}}{\varepsilon_{k}}\omega_k}.
\end{align*}
If we write
$\eta_{k+1}:=\frac{\varepsilon_k}{\varepsilon_{k+1}}\omega_{k+1}, k\ge 1$ and for $k\ge 2,$
$$
  \tilde\beta_{k+1}:=\frac{\delta_k}{\delta_{k+1}}(\alpha_{k+1}
  +\omega_{k}+\varepsilon_{k})\omega_{k+1},\ \tilde\alpha_{k} :=\alpha_k+\omega_{k-1}+\varepsilon_{k-1}
  -\frac{\delta_k}{\delta_{k+1}}\omega_{k+1},
$$
then \begin{align}\label{etk}
  \eta_{k+1}=\frac{\tilde\beta_{k+1}}{\tilde\alpha_{k}+\eta_k}, k\ge2
\end{align}
and
\begin{align*}
  \lim_{k\rto} \tilde\beta_k=\beta q,\ \lim_{k\rto}\tilde\alpha_k=\alpha+\omega-q\omega.
\end{align*}
Applying \cite[Theorem 4.13, page 188]{lw}, the limit $\eta:=\lim_{k\rto}\eta_k$ exists.
 Letting $k\rto$ in (\ref{etk}), we get $\eta=\frac{\beta q}{\alpha+\omega-\omega q +\eta}$ whose  solutions are $$\eta=\omega q \text{ or }\eta=-(\alpha+\omega).$$ Since
$ \eta=\lim_{k\rto}\eta_{k+1}=\lim_{k\rto}\frac{\varepsilon_{k}}{\varepsilon_{k+1}}\omega_{k+1} =\omega\lim_{k\rto}\frac{\varepsilon_{k}}{\varepsilon_{k+1}}$ and $\omega>0,$
we have $$\lim_{k\rto}\frac{\varepsilon_{k}}{\varepsilon_{k+1}}=q\text{ or } -\frac{\alpha+\omega}{\omega}.$$ The lemma is proved.\qed

The lemma below ensures the existence of the limit of $ \frac{\delta_k}{\delta_{k+1}}$ as $k\rto,$ so that the fluctuation of $f_k-\omega_k,k\ge1$ can be studied.
\begin{lemma}\label{eq} Suppose condition (B1) and one of (B2)$_a,$ (B2)$_b$ and  (B2)$_c$ hold. Set $\beta_k=d_k/b_k,\alpha_k=a_k/b_k,k\ge1.$  Let $\omega_k:=\frac{\sqrt{\alpha_{k+1}^2+4\beta_{k+1}}-\alpha_{k+1}}{2},k\ge1$ and $\delta_{k}=\beta_k-\omega_{k}(\alpha_k+\omega_{k-1}),k\ge2.$ Then  $\lim_{k\rto}\frac{\delta_k}{\delta_{k+1}}$ exists as an extending number.
  \end{lemma}
\proof By condition (B1), we have $\beta_k\rightarrow d/b=:\beta,$ $\alpha_k\rightarrow  a/b=:\alpha$ and $\omega_k\rightarrow \frac{\sqrt{\alpha^2+4\beta}-\alpha}{2}=:\omega,$ as $k\rto.$

For simplicity, write temporarily $S_k:=\sqrt{\alpha_k^2+4\beta_k}-\alpha_k,k\ge 1.$ Then we have $\lim_{k\rto} S_k=\sqrt{\alpha^2+4\beta}-\alpha.$ Some direct computation yields that
\begin{align}
  \frac{\delta_k}{\delta_{k+1}}&=\frac{\beta_k}{\beta_{k+1}}\times\frac{\frac{\alpha_{k+1}-\alpha_k}{S_k}
  +\frac{(\alpha_k-\alpha_{k+1})(\alpha_k+\alpha_{k+1})+4(\beta_k-\beta_{k+1})}{S_k\z(S_k+S_{k+1}+\alpha_k+\alpha_{k+1}\y)}}
  {\frac{\alpha_{k+2}-\alpha_{k+1}}{S_{k+1}}
  +\frac{(\alpha_{k+1}-\alpha_{k+2})(\alpha_{k+1}+\alpha_{k+2})+4(\beta_{k+1}-\beta_{k+2})}{S_{k+1}\z(S_{k+1}+S_{k+2}+\alpha_{k+1}+\alpha_{k+2}\y)}},\no\\
  &=:\frac{\beta_k}{\beta_{k+1}}\times \frac{D_k}{D_{k+1}}, k\ge2. \label{dd}
\end{align}
Suppose now condition (B2)$_a$ holds. Then for some $k_0>0,$ $\alpha_k=\alpha_{k+1},\beta_k\ne \beta_{k+1}, \forall k\ge k_0$ and  $$\lim_{k\rto}\frac{\beta_k-\beta_{k+1}}{\beta_{k+1}-\beta_{k+2}}=\lim_{k\rto}\frac{d_{k}/b_{k}-d_{k+1}/b_{k+1}}{d_{k+1}/b_{k+1}-d_{k+2}/b_{k+2}}$$ exists as an extending number. Therefore
\begin{align*}
  \lim_{k\rto}\frac{D_k}{D_{k+1}}&=\lim_{k\rto}\frac{\beta_k-\beta_{k+1}}{\beta_{k+1}-\beta_{k+2}}
  \frac{S_{k+1}\z(S_{k+1}+S_{k+2}+\alpha_{k+1}+\alpha_{k+2}\y)}{S_k\z(S_k+S_{k+1}+\alpha_k+\alpha_{k+1}\y)}\no\\
  &=\lim_{k\rto}\frac{\beta_k-\beta_{k+1}}{\beta_{k+1}-\beta_{k+2}}
\end{align*}
exists as an extending number. Consequently, by (\ref{dd}),  $\lim_{k\rto}\frac{\delta_k}{\delta_{k+1}}$ exists as an extending number.

If condition (B2)$_b$ holds, a similar argument also yields that  $\lim_{k\rto}\frac{\delta_k}{\delta_{k+1}}$ exists.

Suppose condition (B2)$_c$ holds. Then $\exists k_0>0$ such that $\alpha_k\ne \alpha_{k+1}, \beta_k\ne \beta_{k+1},\forall k\ge k_0$ and
\begin{align}
  \tau&:=\lim_{k\rto}\frac{\beta_k-\beta_{k+1}}{\alpha_{k}-\alpha_{k+1}}=\lim_{k\rto}\frac{d_{k}/b_k-d_{k+1}/b_{k+1}}{a_{k}/b_k-a_{k+1}/b_{k+1}}\ne \frac{-a+ \sqrt{a^2+4bd}}{2b}\label{tau} \end{align}
  exists as an extending number.
Furthermore if $\tau$ is finite, by assumption, we also have that \begin{align}\label{lma}  \lim_{k\rto}\frac{\alpha_k-\alpha_{k+1}}{\alpha_{k+1}-\alpha_{k+2}}=\lim_{k\rto}\frac{a_{k}/b_{k}-a_{k+1}/b_{k+1}}{a_{k+1}/b_{k+1}-a_{k+2}/b_{k+2}} \end{align} exists as an extending number.
In this situation, taking (\ref{tau}) into account, we have \begin{align}
  \lim_{k\rto}&\frac{(\alpha_k+\alpha_{k+1})+4(\beta_k-\beta_{k+1})/(\alpha_k-\alpha_{k+1})}{\z(S_k+S_{k+1}+\alpha_k+\alpha_{k+1}\y)}-1\no\\
&=\frac{\alpha+2\tau}{\sqrt{\alpha^2+4\beta}}-1\ne 0\label{ne0}
\end{align}
since $\tau\ne \frac{-a+ \sqrt{a^2+4bd}}{2b}.$ Therefore, by (\ref{tau}), \eqref{lma} and (\ref{ne0}),
\begin{align}\label{dfd}
  &\lim_{k\rto}\frac{D_k}{D_{k+1}}\no\\
  &\quad=\lim_{k\rto}\frac{\alpha_k-\alpha_{k+1}}{\alpha_{k+1}-\alpha_{k+2}}\frac{S_{k+1}}{S_k}
  \frac{\frac{(\alpha_k+\alpha_{k+1})+4(\beta_k-\beta_{k+1})/(\alpha_k-\alpha_{k+1})}{\z(S_k+S_{k+1}+\alpha_k+\alpha_{k+1}\y)}-1}
  {\frac{(\alpha_{k+1}+\alpha_{k+2})+4(\beta_{k+1}-\beta_{k+2})/(\alpha_{k+1}-\alpha_{k+2})}{\z(S_{k+1}+S_{k+2}+\alpha_{k+1}+\alpha_{k+2}\y)}-1}\no\\
  &\quad =\lim_{k\rto}\frac{\alpha_k-\alpha_{k+1}}{\alpha_{k+1}-\alpha_{k+2}}\no
\end{align}
 exists as an extending number.
 Consequently, we conclude from (\ref{dd}) that $$\lim_{k\rto}\frac{\delta_k}{\delta_{k+1}}=\lim_{k\rto}\frac{D_k}{D_{k+1}}=\lim_{k\rto}\frac{\alpha_k-\alpha_{k+1}}{\alpha_{k+1}-\alpha_{k+2}}$$
exists as an extending number.

Finally,  if $\tau=\infty,$  then we must have $\lim_{k\rto}\frac{\alpha_{k}-\alpha_{k+1}}{\beta_k-\beta_{k+1}}=0$
and by assumption,  $$\lim_{k\rto}\frac{\beta_k-\beta_{k+1}}{\beta_{k+1}-\beta_{k+2}}=\lim_{k\rto}\frac{d_{k}/b_{k}-d_{k+1}/b_{k+1}}{d_{k+1}/b_{k+1}-d_{k+2}/b_{k+2}}$$ exists as an extending number.
Thus, \begin{align*}
  \lim_{k\rto}\frac{D_k}{D_{k+1}}&=\lim_{k\rto}\frac{\beta_k-\beta_{k+1}}{\beta_{k+1}-\beta_{k+2}}
  \frac{S_{k+1}\z(S_{k+1}+S_{k+2}+\alpha_{k+1}+\alpha_{k+2}\y)}{S_k\z(S_k+S_{k+1}+\alpha_k+\alpha_{k+1}\y)}\no\\
  &=\lim_{k\rto}\frac{\beta_k-\beta_{k+1}}{\beta_{k+1}-\beta_{k+2}}
\end{align*} exists as an extending number.
As a result,  by (\ref{dd}), $\lim_{k\rto}\frac{\delta_k}{\delta_{k+1}}$ exists as an extending number. The lemma is proved. \qed

\subsection{Proof of Theorem \ref{amij}}
We are now ready to prove Theorem \ref{amij}. To begin with, for $k\ge m\ge1$  we write
$x_{k,m}:=\frac{\mathbf e_1A_k\cdots A_m\mathbf e_1^t}{\varrho(A_k)\cdots\varrho(A_m)}$ for simplicity. In view of (\ref{fm}) and (\ref{gm}), to prove \eqref{rac}, it suffices to show that $x_{k,m}\rightarrow c(m)$ as $k\rto$ for some $0<c(m)<\infty.$

It follows from Lemma \ref{lwb} that there exist some constants $0<c_3<c_4<\infty$ independent of $k$ and $m$ such that
$ c_3\le x_{k,m}\le c_4,\forall k\ge m\ge1.$
Therefore, if \eqref{rac} is true,  we must have  $c_3\le c(m)\le c_4,m\ge1.$ Thus the second part of Theorem \ref{amij} holds.

Next, prove only $x_{k,1}\rightarrow c $ as $k\rto$ for some $0<c<0,$  since for $m\ge2$ the convergence of $x_{k,m}$ as $k\rto$ can be proved similarly. In the remainder of the paper, we write $x_{k,1}$ simply as $x_k$ and keep in mind that
\begin{equation}\label{xul}
  x_{k}=\frac{\mathbf e_1A_k\cdots A_1\mathbf e_1^t}{\varrho(A_k)\cdots\varrho(A_1)} \text{ and } \forall k\ge1, c_3\le x_k\le c_4.
\end{equation}

Let $f_k,k\ge1$ be the one in (\ref{ctf}) with $\beta_k=d_k/b_k,\alpha_k=a_k/b_k,k\ge1$ and set $\omega_k:=\frac{\sqrt{a_{k+1}^2+4b_{k+1}d_{k+1}}-a_{k+1}}{2b_{k+1}},k\ge 1.$
 Then
 $$\beta_k\rightarrow d/b>0, \alpha_k\rightarrow a/b>0, \omega_k \rightarrow \frac{\sqrt{a^2+4bd}-a}{2b}=:\omega>0$$ as $k\rto.$ Let $\delta_k,k\ge1$ and $\varepsilon_k,k\ge2$ be those in \eqref{dlep}.

By (\ref{fm}), (\ref{cfa}) and (\ref{ff}), we have $\lim_{k\rto}f_k=\frac{\sqrt{a^2+4bd}-a}{2b}=\omega,$ implying that $\lim_{k\rto}\varepsilon_k=0.$

 An application of  Lemma \ref{eq} yields that the limit $\lim_{k\rto}\frac{\delta_k}{\delta_{k+1}}$ exists as an extending number. Assume $$\lim_{k\rto}\frac{\delta_k}{\delta_{k+1}}=q.$$ Since $\lim_{k\rto}\delta_k=0,$ we must have $|q|\ge1.$ Applying  Lemma \ref{de}, we have
\begin{align}\label{elk}
  \lim_{k\rto}\frac{\varepsilon_k}{\varepsilon_{k+1}}=q\text{ or }-\frac{a/b+\omega}{\omega}.
\end{align}

{\it Case 1: Suppose $|q|>1.$} Then $\lambda_0:=\min\{|q|,\frac{a/b+\omega}{\omega}\}>1.$ Fix some $1<\lambda<\lambda_0.$ By (\ref{elk}), there exists some $k_0>0$ such that $\z|\frac{\varepsilon_k}{\varepsilon_{k+1}}\y|\ge \lambda, $ for all $k>k_0.$ Hence
\begin{equation}\label{foa}
  \sum_{k=2}^{\infty}|\varepsilon_k|=\sum_{k=2}^\infty|f_k-\omega_k| <\infty.
\end{equation}
Taking \eqref{rhk} and (\ref{cfa}) into account, for $k\ge 1$ we have
\begin{align}\label{xf}
  x_{k+1}-x_k&=\frac{a_{k+1}-\varrho(A_{k+1})+b_{k+1}\frac{\mathbf e_2A_k\cdots A_1\mathbf e_1^t}{\mathbf e_1A_k\cdots A_1\mathbf e_1^t}}{\frac{\varrho(A_{k+1})\cdots\varrho(A_1)}{\mathbf e_1A_k\cdots A_1\mathbf e_1^t}}\no\\
  &=\z(\frac{\varrho(A_{k+1})\cdots\varrho(A_1)}{\mathbf e_1A_k\cdots A_1\mathbf e_1^t} \y)^{-1}(a_{k+1}-\varrho(A_{k+1})+b_{k+1}f_k)\no\\
  &=(\varrho(A_{k+1})x_k^{-1})^{-1}b_{k+1}(f_k-\omega_k).
\end{align}
Since $\varrho(A_{k}),k\ge 1$ is uniformly bounded away from $0$ and $\infty,$ then by (\ref{xul}), we have for some constant $0<c_5<\infty,$
\begin{align}\label{xda}
  |x_{k+1}-x_k|\le c_5 |f_k-\omega_k|, \forall k\ge1.
\end{align}
Taking (\ref{xul}), (\ref{foa}) and (\ref{xda}) together, we conclude that for some constant $0<c<\infty,$ $\lim_{k\rto} x_k=c.$

{\it Case 2: Suppose $q=1$ and $\lim_{k\rto}\frac{\varepsilon_k}{\varepsilon_{k+1}}=-\frac{a/b+\omega}{\omega}.$ } Since $\frac{a/b+\omega}{\omega}>1,$ the proof goes exactly the same as Case 1.

{\it Case 3: Suppose $q=1$ and $\lim_{k\rto}\frac{\varepsilon_k}{\varepsilon_{k+1}}=q.$} Then there exists some number $k_1>0$ such that $\varepsilon_k=f_k-\omega_k, k\ge k_1$ are all strictly positive or strictly negative, and consequently
\begin{align}\label{skof}\frac{a_{k+1}+b_{k+1}f_k}{a_{k+1}+b_{k+1}\omega_k}<1(\text{or }>1), \text{ for all }k\ge k_1.\end{align}
But \begin{align}\label{ratio}
  \frac{x_{k+1}}{x_k}&=\frac{1}{\varrho(A_{k+1})}\z(a_{k+1}+b_{k+1}\frac{\mathbf e_2A_k\cdots A_1\mathbf e_1^t}{\varrho(A_k)\cdots\varrho(A_1)}\y)\no\\
  &=\frac{1}{\varrho(A_{k+1})}\z(a_{k+1}+b_{k+1}f_k\y)=\frac{a_{k+1}+b_{k+1}f_k}{a_{k+1}+b_{k+1}\omega_k}.
\end{align}
Thus, by (\ref{skof}), $\frac{x_{k+1}}{x_k}<1(\text{or }>1)$ for all $k\ge k_1,$ that is, $x_k,k\ge k_1$ is monotone. As a consequence, it follows from (\ref{xul}) that for some constant $0<c<\infty$ $\lim_{k\rto}x_k=c.$

{\it Case 4. Suppose that $q=-1$ and $\lim_{k\rto}\frac{\varepsilon_k}{\varepsilon_{k+1}}=-\frac{a/b+\omega}{\omega}.$} In this case, the proof is the same as Case 2.

{\it Case 5. Suppose that $q=-1$ and $\lim_{k\rto}\frac{\varepsilon_k}{\varepsilon_{k+1}}=-1.$}
Combining (\ref{xf}) with \eqref{ratio}, we have $$ \frac{x_{k+1}-x_k}{x_{k}-x_{k-1}}=\frac{\varrho(A_k)}{\varrho(A_{k+1})}
\frac{a_k+b_k\omega_{k-1}}{a_k+b_kf_{k-1}}\frac{b_{k+1}}{b_k}\frac{\varepsilon_{k}}{\varepsilon_{k-1}}\rightarrow -1,$$ as $k\rto.$
 So there exists some number $k_2>0$ such that
 $$ \frac{x_{k+1}-x_k}{x_{k}-x_{k-1}}<0\ \text{ for all } k> k_2.$$
Since $\varepsilon_k= f_k-\omega_k\rightarrow 0$ as $k\rto,$ then by (\ref{xf}), we have $x_{k+1}-x_k\rightarrow0$ as $k\rto.$
We thus come to the conclusion that $x_{k+1}-x_k$ converges to $0$ in an alternating manner as $k\rto.$ Therefore, $$c:=\lim_{k\rto}x_k=x_1+\sum_{k=1}^{\infty} (x_{k+1}-x_k)$$ exists and by (\ref{xul}), we must have $0<c<\infty.$ Theorem \ref{amij} is proved.\qed

\section{Maximum of (2,1) random walk }\label{mrw}

\subsection{Escape probability and the distribution of $M$}\label{sba}
We consider firstly the escape probability of (2,1) random walk $Y$ from certain interval.  For $0<m\le k\le n,$  let
\begin{equation*}
  P_k(m,n,-)=P(Y\text{ hits }[0,m]\text{ before it hits }[n,\infty)|Y_0=k).
\end{equation*}
The above escape probabilities can be written as functionals of product of nonnegative matrices.
\begin{lemma}\label{esy} For $0<m\le k\le n,$ we have
  \begin{align}\label{epp}P_k(m,n,-)=\frac{\sum_{s=k}^{n-1}\mathbf e_1 N_s\cdots N_{m+1}\mathbf e_1^t}{1+\sum_{s=m+1}^{n-1}\mathbf e_1 N_s\cdots N_{m+1}\mathbf e_1^t}.\end{align}
\end{lemma}
\proof The proof of the lemma is very standard and can be find in \cite{br}. Here we sketch its proof for convenience of the reader.  For $j\in\{m,m-1\},$ let $$P_k^j(m,n,-):=P(Y\text{ hits }[0,m]\text{ at } j\text{ before it hits }[n,\infty)|Y_0=k)$$ which we will write as $P_k^j$ for simplicity. Then by Markov property,
$$P_k^j=q_kP^j_{k+1}+p_kP_{k-2}^j,\ m+1\le k\le n-1,$$ which leads to
 \begin{align}
  P_{k+1}^j-P_k^j=\frac{p_k}{q_k}(P_{k}^j-P_{k-1}^j)+\frac{p_k}{q_k}(P_{k-1}^j-P_{k-2}^j)\label{dif}
 \end{align}
 with the boundary condition \begin{align}P_m^m=P_{m-1}^{m-1}=1, P_m^{m-1}=P_{m-1}^m=P_n^{m}=P_n^{m-1}=0.\label{il}\end{align}
Set $V_k^j=\z(\begin{array}{c}
                P_{k+1}^j-P_k^j \\
                P_{k}^j-P_{k-1}^j
              \end{array}
\y). $ Clearly, \begin{align}
  V_m^m=\z(\begin{array}{c}
                P_{m+1}^m-1 \\
                1
              \end{array}\y),
V_m^{m-1}=\z(\begin{array}{c}
                P_{m+1}^{m-1}\\
                -1
              \end{array}\y).\label{bdc}
\end{align}
By (\ref{dif}), we get $
  V_k^j=\z(\begin{array}{cc}
                          p_k/q_k & p_k/q_k \\
                           1& 0
                        \end{array}
\y)V_{k-1}^j=N_kV_{k-1}^j,
$
implying that
\begin{align}\label{vnk}
  V_k^j=N_kN_{k-1}\cdots N_{m+1}V_{m}^j,\ m+1\le k\le n-1.
\end{align}
Taking (\ref{il}) and the fact $P_k^m(m,n,-)+P_k^{m-1}(m,n,-)=P_{k}(m,n,-)$ into account, solving (\ref{vnk}) with initial condition (\ref{bdc}), we get  (\ref{epp}). \qed

With Lemma \ref{esy} in hands, we next prove Proposition \ref{dm}.

\noindent{\it Proof of Proposition \ref{dm}. } Note that the event $\{M=n,D<\infty\}$ occurs if and only if the following two events occur successively: i) starting from $2,$ the walk hits $n$ before the set $\{0,1\}$ and then ii) restarting from $n,$  it hits the set $\{0,1\}$ before $n+1.$ Thus, it follows from Markov property that $$P(M=n,D<\infty)=(1-P_2(1,n,-))P_n(1,n+1,-).$$
Consequently, an application of Lemma \ref{esy}  finishes the proof of (\ref{pmn}).  \qed

\subsection{Asymptotics of $\prod_{i=2}^k\varrho(N_i)$ and the distribution of $M$}\label{sbb}
Recall that $
N_k:=\z(\begin{array}{cc}
          \theta_k & \theta_k \\
          1 & 0
        \end{array}
\y)\text{ with }\theta_k:=\frac{p_k}{q_k}.$ Our aim is to prove Theorem \ref{mxt}. For this purpose, since by Proposition \ref{dm}, the distribution of $M$ is written in terms of $\mathbf e_1 N_k\cdots N_{2}\mathbf e_1^t, k\ge2,$ which is hard to estimate directly.
 But by Theorem \ref{amij}, it is sufficient to work with $\varrho(N_k)\cdots \varrho(N_1).$

  Note that $\varrho(N_k)=\z(\theta_k+\sqrt{\theta_k^2+4\theta_k}\y)/2.$ If $q_i=\frac{2}{3}\pm r_i,i\ge 2,$ then by Taylor enpension of  $\varrho(N_k)$ at $0,$ we get  \begin{equation}\varrho(N_k)=1\mp 3r_k+O(r_k^2)\text{ as }k\rto.\label{rs}\end{equation}
 The proposition below yields the asymptotics of $\varrho(N_k)\cdots \varrho(N_1).$

 \begin{proposition}\label{pra} Fix $K= 1,2,3,...$ and $B\in \mathbb R.$ Let $r_i$ be the one in \eqref{ri}. Suppose that $\sigma_i,i\ge 2$ is a sequence of numbers such that
 \begin{align}\label{sia}
   \sigma_i=1\pm 3r_i+O(r_i^2)\text{ as }i\rto.
 \end{align}
Then we have
 \begin{align}\label{ra}
  &\sigma_2\cdots \sigma_n\sim  c\z(n\log n\cdots \log_{K-2}n(\log_{K-1} n)^B\y)^{\pm 1}\end{align}
  and
  \begin{align}\label{sto}
  &\frac{\sigma_2\cdots \sigma_n}{\sum_{i=1}^{n}\sigma_2\cdots \sigma_i}\rightarrow 0
\end{align}
as $n\rto.$
\end{proposition}
 \proof
 Assume $\sigma_i=1\pm 3r_i+O(r_i^2)\text{ as }i\rto.$ For a rigorous proof of (\ref{ra}), we refer the reader to \cite[Lemma 2]{wh}. Here, we only sketch its proof. For some $M>0$ sufficiently large, we have
 \begin{align*}
  \sigma_2\cdots \sigma_n&=\exp\z\{\sum_{i=1}^n\log\sigma_i\y\}\\
  &= \exp\z\{\sum_{i=1}^n\log(1\pm 3r_i+O(r_i^2))\y\}\sim c\exp\z\{\pm\sum_{i=1}^n3r_i\y\}\\
  &=c\exp\z\{\pm\sum_{i=1}^n\Lambda(K,i,B)\y\}\sim c \exp\z\{\pm\int_{M}^n\Lambda(K,i,B)\y\}\\
  & \sim c\z(n\log n...\log_{K-2}n(\log_{K-1}n)^B\y)^{\pm1}, \ \mbox{as} \ n\ \rto,
 \end{align*}
 which leads to (\ref{ra}).

Finally,  \eqref{sto} is a direct consequence of \eqref{ra}. \qed

 For the product $N_k\cdots N_2,k\ge2,$   requirements of Theorem \ref{amij} are fulfilled by the following lemma.
\begin{lemma}\label{dfr}{\rm (i)} We have $\lim\limits_{n\rightarrow\infty}\frac{r_n-r_{n+1}}{n^2}=1/3$ and thus $\sum_{k=2}^\infty|\theta_{k+1}-\theta_k|<\infty$ whenever $q_i=2/3\pm r_i, i\ge 2.$
{\rm (ii)} for $k\ge i_0,$ we have $\frac{1}{\theta_k}\neq\frac{1}{\theta_{k+1}}$ and
$$\lim_{k\rto}\frac{\theta_{k+1}-\theta_k}{\theta_{k+2}-\theta_{k+1}}=1,$$ no mater $q_i=2/3+r_i$ or $q_i=2/3-r_i,i\ge 2.$

\end{lemma}
\proof  Fix $K=1,2,3,...,$ $B\in\mathbb R$ and $n\ge i_0.$ For $i\ge k\ge 0,$ set
\begin{align*}
  &\Gamma_i=\frac{1}{\prod_{j=0}^i\log_jn}-\frac{1}{\prod_{j=0}^i\log_j(n+1)},\no\\
  & I_k= \frac{1}{\prod_{j=0}^{k-1}\log_j(n+1)\prod_{j=k}^i\log_jn}
  -\frac{1}{\prod_{j=0}^{k}\log_j(n+1)\prod_{j=k+1}^i\log_jn}.
\end{align*}
Here and throughout the paper, we use the convention that empty product equals identity.

Obviously, we have
 \begin{align*}
 r_n-r_{n+1}=\frac{1}{3}\z(\Gamma_0+\cdots+\Gamma_{K-2}+B\Gamma_{K-1}\y),
 \end{align*}
  and for $0\le i\le K-1$,
\begin{align*}
\Gamma_i=I_0+I_1+\cdots+I_i.
\end{align*}
Now, fix $1\le i\le K-1.$ Clearly we have
$$I_0=\frac{1}{n(n+1)\log n\times\cdots\times\log_{i}n}=o\z(\frac{1}{n^2}\y),\ n\rightarrow\infty,$$
and for $1\le k\le i,$
by the mean value theorem,
\begin{align*}
I_k
&=\frac{1}{\prod_{j=0}^{k-1} \log_j(n+1)\prod_{j=k+1}^{i}\log_{j}n}\Big(\frac{1}{\log_{k}n}-\frac{1}{\log_{k}(n+1)}\Big)\\
&=\frac{1}{\prod_{j=0}^{k-1} \log_j(n+1)\prod_{j=k+1}^{i}\log_{j}n}\frac{1}{\prod_{j=0}^{k-1}\log_j\theta_n\log_k^2\theta_n},
\end{align*}
with some $ \theta_n\in(n,n+1).$
It is easy to see that $$\frac{1}{\prod_{j=0}^{k-1}\log_j\theta_n\log_k^2\theta_n}
\sim\frac{1}{\prod_{j=0}^{k-1}\log_jn\log_k^2n}=o\z(\frac{1}{n}\y),\text{ as }n\rightarrow\infty.$$
Thus, $I_k=o(\frac{1}{n^2}), \ \forall 1\le k\le i.$ Consequently, we get
 $$\forall 1\le i\le K-1,\ \Gamma_i=o\z(\frac{1}{n^2}\y),\text{ as }n\rightarrow\infty.$$
Notice  that $\Gamma_0=\frac{1}{n}-\frac{1}{n+1}\sim \frac{1}{n^2}$ as $n\rto.$
Thus we conclude that $$r_n-r_{n+1}\sim \frac1{3n^2},\text{ as } n\rightarrow\infty.$$

On the other hand,  when $q_i=2/3\pm r_i,i\ge2,$ $$\sum_{k=2}^\infty|\theta_{k+1}-\theta_k|\le c\sum_{k=2}^\infty|r_{k+1}-r_k|\le c\sum_{k=2}^{\infty}\frac{1}{k^2}<\infty.$$

To prove the second part, noting that $r_k\ne r_{k+1}, \forall k\ge i_0,$  thus $\theta_{k+1}^{-1}\ne \theta_{k}^{-1}, \forall k\ge i_0.$
Note also that if $q_i=2/3\pm r_i,i\ge 2,$ then $\theta_k=\frac{p_k}{q_k}=\frac{1}{2}\pm \frac{9}{4}r_i+\frac{27}{4}r_k^2+o(r_k^2)$ as $k\rto.$
Since $r_n-r_{n+1}\sim \frac1{3n^2}$ as $ n\rightarrow\infty $ by the first part, then we have  $\lim_{k\rto}\frac{\theta_{k+1}-\theta_k}{\theta_{k+2}-\theta_{k+1}}=1.$
The lemma is proved. \qed

\subsection{Proof of Theorem \ref{mxt}}\label{pta}
To conclude Section \ref{mrw}, we give the proof of Theorem \ref{mxt}.
 Recall that by Proposition \ref{dm}, we have
\begin{align*}
    P(M=n,D<\infty)
    =\frac{1}{1+\sum_{s=2}^{n-1}\mathbf e_1 N_s\cdots N_{2}\mathbf e_1^t}\frac{\mathbf e_1 N_n\cdots N_{2}\mathbf e_1^t}{1+\sum_{s=2}^{n}\mathbf e_1 N_s\cdots N_{2}\mathbf e_1^t}.
  \end{align*}
Consider the product of matrices $N_k\cdots N_2.$ Lemma \ref{dfr} ensures that conditions (B1) and (B2)$_a$ are fulfilled. Thus, applying Theorem \ref{amij}, we get that $\mathbf e_1N_k\cdots N_2\mathbf e_1^t\sim c\varrho(N_k)\cdots \varrho(N_2),$ as $k\rto.$

If $q_i=\frac{2}{3}\pm r_i,i\ge 2,$ then by \eqref{rs},  $\varrho(N_k)=1\mp 3r_k+O(r_k^2)\text{ as }k\rto.$ Applying Proposition \ref{pra}, we get
$$\varrho(N_2)\cdots\varrho(N_n)\sim c\z(n\log n\cdots \log_{K-2}n(\log_{K-1} n)^B\y)^{\mp 1}$$ as $n\rto.$
Thus $$\mathbf e_1N_n\cdots N_2\mathbf e_1^t\sim c\z(n\log n\cdots \log_{K-2}n(\log_{K-1} n)^B\y)^{\mp 1},$$ as $n\rto.$
Consequently, the proof of Theorem \ref{mxt} is  an almost verbatim repetition of \cite[Theorem 1]{wh}. We do not repeat it here. \qed

\section{Maximum of (1,2) random walk}\label{mt}
In this section, we consider (1,2) random walk $Y',$ whose escape probabilities can be written in terms of the tails of a continued fraction. By analyzing the escape probabilities and tails of the continued fraction delicately, we can study the asymptotics of $M(Y').$

\subsection{Continued fraction and escape probability}\label{mt1}
Let $\theta_k,k\ge2$ be those defined in \eqref{tn}. Consider the continued fraction
\begin{equation*}\label{cfth}
\mathrm K_{n=2}^\infty(\theta_n^{-1}|1)\equiv\frac{\theta_2^{-1}}{1}\begin{array}{c}
                                \\
                               +
                             \end{array}\frac{\theta_3^{-1}}{1}\begin{array}{c}
                                \\
                               +
                             \end{array}\frac{\theta_4^{-1}}{1}\begin{array}{c}
                                \\
                               +\cdots
                             \end{array}
\end{equation*}
whose $n$-th tail is denoted by
\begin{align}
  f^{(n)}=\mathrm K_{k=n+1}^\infty(\theta_k^{-1}|1)\equiv\frac{\theta_{n+1}^{-1}}{1}\begin{array}{c}
                                \\
                               +
                             \end{array}\frac{\theta_{n+2}^{-1}}{1}\begin{array}{c}
                                \\
                               +
                             \end{array}\frac{\theta_{n+3}^{-1}}{1}\begin{array}{c}
                                \\
                               +\cdots
                             \end{array}\label{tail}
\end{align}
It would be convenient to write
\begin{equation}
  \xi_{n+1}=f^{(n)},\text{ for }n\ge 1.\label{dx}
\end{equation}

The lemma below presents several limit behaviors related to $\xi_n,n\ge2.$
\begin{lemma}\label{xie}
If $p_i=1/3\pm r_i,i\ge 2,$  then we have
  \begin{equation}\label{xir}\xi_{n}=1\mp 3r_{n}+O(r_{n}^2) \text{ as }n\rto,\end{equation}
  and consequently, \begin{align}\label{xa}
\xi_2\cdots \xi_n\sim  c&\z(n\log n\cdots \log_{K-2}n(\log_{K-1} n)^B\y)^{\mp 1},\\
  &\frac{\xi_2\cdots \xi_n}{\sum_{i=1}^{n}\xi_2\cdots \xi_i}\rightarrow 0\label{xto}
\end{align}
as $n\rto.$
\end{lemma}
\proof Since $r_i,i\ge n_0$ is monotone for some $n_0>0$ large enough and by Lemma \ref{dfr}, $(r_n-r_{n+1})\sim 3r_n^2\sim\frac{1}{3n^2}$ as $n\rto,$ then the proof of \eqref{xir} goes step by step as that of \cite[Lemma 1]{w19}.

By \eqref{xir}, the condition of Proposition \ref{pra} is satisfied. Thus, apply Proposition \ref{pra} to the sequence $\xi_i,i\ge 2,$ we obtain \eqref{xa} and \eqref{xto}. \qed

For  $1\le m\le k\le n,$ $j\in \{n,n+1\},$  let
\begin{align*}
&\mathcal P_k^j(m,n,+):=P(Y'\text{ hits }[n,\infty] \text{ at }j\text{ before }[0,m] |Y'_0=k).\\
&\mathcal P_k(m,n,+):=P(Y'\text{ hits }[n,\infty] \text{ before }[0,m] |Y'_0=k)\\
&\mathcal P_k(m,n,-):=1-\mathcal P_{k}(m,n,+).
\end{align*}
Clearly, we have $\mathcal P_k(m,n,+)=\mathcal P_k^n(m,n,+)+\mathcal P_k^{n+1}(m,n,+).$
We have the following estimations of the escape probabilities.
\begin{lemma}\label{ecpb} For any integers $1\le m\le k\le n,$

\begin{equation}\frac{\sum_{i=k}^{n-1}\xi_{m+1}\cdots \xi_i}{1+\sum_{i=m+1}^{n-1}\xi_{m+1}\cdots \xi_i}\le \mathcal P_k(m,n,-)\le \frac{\sum_{i=k}^{n}\xi_{m+1}\cdots \xi_i}{1+\sum_{i=m+1}^{n}\xi_{m+1}\cdots \xi_i}.\label{exi}\end{equation}
\end{lemma}
The lemma can be proved by a space reversal argument of  the proof of \cite[Lemma 1 on page 230]{leta}.

\subsection{Ratio of the hitting probabilities} \label{mt2}

It is easily seen that
\begin{align}P(M=n, D<\infty)=\mathcal P_2^n(1,n,+)\mathcal P_n(1,n+1,-),\ n\ge 2.\no\end{align}
 From \eqref{exi}, we can get
 \begin{align}
   \frac{\xi_{2}\cdots \xi_{n}}{1+\sum_{i=2}^{n}\xi_{2}\cdots \xi_i}\le \mathcal P_n(1,n&+1,-)\le \frac{\xi_{2}\cdots \xi_n+\xi_{2}\cdots \xi_{n+1}}{1+\sum_{i=2}^{n+1}\xi_{2}\cdots \xi_i},\label{esa}\\
   \frac{1}{1+\sum_{i=2}^{n}\xi_{2}\cdots \xi_i}\le \mathcal P_2(1,&n,+)\le \frac{1}{1+\sum_{i=2}^{n-1}\xi_{2}\cdots \xi_i}.\label{esb}
 \end{align}
 But  it is impossible to compute or estimate  $\mathcal P_2^n(1,n,+)$ directly from \eqref{exi}. The following lemma shows that the ratio of $\mathcal P_2^n(1,n,+)$ over $\mathcal P_2^{n+1}(1,n,+)$ has a limit.

\begin{lemma}\label{rnn}  Suppose that $p_i=\frac{1}{3}\pm r_i,i\ge2.$ Then we have
\begin{align}
  \lim_{n\rto}\frac{\mathcal P_2^{n}(1,n,+)}{\mathcal P_2^{n+1}(1,n,+)}=2.\label{rpk}
\end{align}
  \end{lemma}
\proof  Fix $n\ge3.$ Let $$E_n=\{Y'\text{ hits }[n,\infty)\text{ before it hits }[0,1]\}.$$
Define a measure $\tilde P$ by $$\tilde P(\cdot)=P(\cdot|E_n).$$
Let $$T_n:=\inf\{k\ge 0: Y_k'\in[n,\infty)\},n\ge3.$$ Then $T_n<\infty$ almost surely.

{\it Step 1.} We claim that $Y'$ is a Markov chain under $\tilde P$ with transition probabilities
\begin{align}
  &\tilde P(Y_{k+1}'=4|Y_k'=2,k<T_n)=1,\label{t24}\\
  &\tilde P(Y_{k+1}'=i+2|Y_k'=i,k<T_n)=p_i\frac{\mathcal P_{i+2}(1,n,+)}{\mathcal P_i(1,n,+)}=:\tilde p_i,\label{t2}\\
  &\tilde P(Y_{k+1}'=i-1|Y_k'=i,k<T_n)=1-\tilde p_i=:\tilde{q}_i,3\le i\le n-1.\no
\end{align}
We first show the Markov property.
Indeed, for $y_j\in [2,n+1], 1\le j< T_n,$ \begin{align}
  \tilde P(Y_{k+1}'&=y_{k+1}|Y_j'=y_j,1\le j\le k<T_n)\no\\
  &= \frac{ P(Y_{k+1}'=y_{k+1},Y_j'=y_j,1\le j\le k<T_n,E_n)}{ P(Y_j'=y_j,1\le j\le k<T_n,E_n)}\no\\
  &=\frac{P(Y_{k+1}'=y_{k+1},E_n |Y_j'=y_j,1\le j\le k<T_n)}{P(E_n|Y_j'=y_j,1\le j\le k<T_n)}\no\\
  &=\frac{P(Y_{k+1}'=y_{k+1},E_n |Y_k'=y_k,k<T_n)}{P(E_n|Y_k'=y_k,k<T_n)}\no\\
  &=\frac{P(Y_{k+1}'=y_{k+1},Y_k'=y_k,k<T_n|E_n)}{P(Y_k'=y_k,k<T_n|E_n)}\no\\
  &=\tilde P(Y_{k+1}'=y_{k+1}|Y_k'=y_k,k<T_n)\no
\end{align}
where the third equality follows from the fact that $Y'$ is a Markov chain under the measure $P.$
Therefore, $Y'$ is also a Markov chain under the measure $\tilde P.$

Next we compute the transition probabilities. Since \eqref{t24} is trivial, we need only to prove \eqref{t2}. It is easy to see that
\begin{align}
  \tilde p_i&=\tilde P(Y_{k+1}'=i+2|Y_k'=i, k<T_n)\no\\
  &=\frac{P(Y_{k+1}'=i+2,Y_k'=i, k<T_n)P(E_n|Y_{k+1}'=i+2,Y_k'=i, k<T_n)}{P(Y_k'=i, k<T_n)P(E_n|Y_k'=i, k<T_n)}\no\\
  &=P(Y_{k+1}'=i+2|Y_k'=i, k<T_n)\frac{P(E_n|Y_{k+1}'=i+2, k<T_n)}{P(E_n|Y_k'=i, k<T_n)}\no\\
  &=p_i\frac{\mathcal P_{i+2}(1,n,+)}{\mathcal P_i(1,n,+)},\no
\end{align}
which proves \eqref{t2}.

{\it Step 2.} We show that $\lim_{i\rto}\tilde p_i=1/3.$ Indeed, by \eqref{exi}, we have
\begin{align}
  &\frac{1+\sum_{j=2}^{i+1}\xi_2\cdots \xi_i}{1+\sum_{j=2}^{i-1}\xi_2\cdots \xi_i}\frac{1+\sum_{j=2}^{n-1}\xi_2\cdots \xi_i}{1+\sum_{j=2}^{n}\xi_2\cdots \xi_i}\no\\
  &\quad\quad\le \frac{\mathcal P_{i+2}(1,n,+)}{\mathcal P_i(1,n,+)}\le \frac{1+\sum_{j=2}^{i+1}\xi_2\cdots \xi_i}{1+\sum_{j=2}^{i-1}\xi_2\cdots \xi_i}\frac{1+\sum_{j=2}^{n}\xi_2\cdots \xi_i}{1+\sum_{j=2}^{n-1}\xi_2\cdots \xi_i}.\label{lup}
\end{align}
Thus, since $i<n,$ it follows from \eqref{xto} that both the leftmost-hand side and the rightmost-hand side of  \eqref{lup} converge to $1$ as $i\rto.$ Therefore, $\lim_{i\rto}\tilde p_i=\lim_{i\rto}p_i=1/3.$

{\it Step 3.} For $2\le k< n,$ define
\begin{align}
  \eta_k(1)=\tilde P(Y' \text{ hits } [k+1,\infty) \text{ at }k+1|Y_0'=k),\no\\
  \eta_k(2)=\tilde P(Y' \text{ hits } [k+1,\infty) \text{ at }k+2|Y_0'=k).\no
\end{align}
Then $\eta_k(1)+\eta_k(2)=1.$
We need to show that
\begin{align}
  \lim_{k\rto}\eta_k(2)=\frac{1}{2}.\label{ek}
\end{align}
To this end, using Markov property, we have
$$\eta_k(2)=\tilde p_k+\tilde q_k\eta_{k-1}(1)\eta_k(2)$$ which leads to
\begin{align}
  \frac{\tilde q_{k+1}}{\tilde p_{k+1}}\eta_{k}(2)=\frac{\frac{\tilde q_{k+1}}{\tilde p_{k+1}}}{1+\frac{\tilde q_k}{\tilde p_k}\eta_{k-1}(2)}.\label{er}
\end{align}
Write $\Gamma_k:=\frac{\tilde q_{k+1}}{\tilde p_{k+1}}\eta_{k}(2)$ and $\beta_k:=\frac{\tilde q_k}{\tilde p_k}.$
Since $\eta_2(2)=1,$ iterating \eqref{er}, we get
\begin{align}
  \Gamma_k=\frac{\beta_{k+1}}{1}\begin{array}{c}
                                \\
                               +
                             \end{array}\frac{\beta_{k}}{1}\begin{array}{c}
                                \\
                               +\cdots  +
                             \end{array}\frac{\beta_{3}}{1}\begin{array}{c}
                             \end{array}, k\ge 3.\no
\end{align}
Since by Step 2, $\lim_{k\rto}\beta_k=2,$ applying again the convergence of limit periodic continued fractions(see \cite[Theorem 4.13, page 188]{lw}), we have that  $\lim_{k\rto}\Gamma_k$ exists. So $\eta:=\lim_{k\rto}\eta_{k}(2)$ exists. Let $k\rto$ in \eqref{er}, we get
$$\eta=\frac{1}{1+2\eta}$$ whose positive solution is $\eta=1/2.$ Thus \eqref{ek} is proved.

{\it Step 4.} For $2\le k\le n,$ let
\begin{align}
 h_k(1)=\tilde{P}(T_k=k|Y_0'=2), h_k(2)=\tilde{P}(T_k=k+1|Y_0'=2).\no
\end{align}
Clearly, we have $h_k(1)+h_k(2)=1.$
We claim that $\lim_{k\rto}h_k(2)=\frac{1}{3}.$

In fact, using Markov property, for $k\ge 3,$ we have
\begin{equation}
  h_k(2)=h_{k-1}(1)\eta_{k-1}(2)=\eta_{k-1}(2)-\eta_{k-1}(2)h_{k-1}(2).\label{it}
\end{equation}
Iterating \eqref{it} and using that fact $h_2(2)=0,$ we obtain
\begin{align}
  h_k(2)=\sum_{i=3}^{k-1}(-1)^{k-i+1}\eta_i\eta_{i+1}\cdots\eta_{k-1}.\label{hk}
\end{align}
Since by \eqref{ek}, $\lim_{i\rto} \eta_i(2)=1/2,$ then for any $1/2>\epsilon>0,$ there exists a number $k_3>0$ such that $1/2-\epsilon< \eta_i<1/2+\epsilon, \forall i\ge k_3.$ It follows from \eqref{hk} that
\begin{align}
  h_k(2)=\sum_{i=3}^{k_3}(-1)^{k-i+1}\eta_i\eta_{i+1}\cdots\eta_{k-1}
  +\sum_{i=k_3+1}^{k-1}(-1)^{k-i+1}\eta_i\eta_{i+1}\cdots\eta_{k-1}.\label{ht}
\end{align}
It is easy to see that
\begin{align}
  \lim_{k\rto}&\z|\sum_{i=3}^{k_3}(-1)^{k-i+1}\eta_i\eta_{i+1}\cdots\eta_{k-1}\y|\le \lim_{k\rto}\sum_{i=3}^{k_3}\eta_i\eta_{i+1}\cdots\eta_{k-1}\no\\
&\le \lim_{k\rto}(k_3-2)(1/2+\epsilon)^{k-k_3-1}=0\no
\end{align}
and
\begin{align}
  &\frac{1/2-\epsilon}{1-\z(1/2-\epsilon\y)^2}- \frac{\z(1/2+\epsilon\y)^2}{1-\z(1/2+\epsilon\y)^2}\no\\
  &\quad\quad\le\varliminf_{k\rto}\sum_{i=k_3+1}^{k-1}(-1)^{k-i+1}\eta_i\eta_{i+1}\cdots\eta_{k-1}\no\\
  &\quad\quad\le\varlimsup_{k\rto}\sum_{i=k_3+1}^{k-1}(-1)^{k-i+1}\eta_i\eta_{i+1}\cdots\eta_{k-1}\no\\
  &\le  \frac{1/2+\epsilon}{1-\z(1/2+\epsilon\y)^2}- \frac{\z(1/2-\epsilon\y)^2}{1-\z(1/2-\epsilon\y)^2}.\no
  \end{align}
  Since $\epsilon$ is arbitrary, the second term  in the righthand side of \eqref{ht} converges to $1/3.$
  Consequently, $\lim_{k\rto}h_k(2)=1/3.$

{\it Step 5.} At last, let us prove $\lim_{n\rto}\frac{\mathcal P_2^{n}(1,n,+)}{\mathcal P_2^{n+1}(1,n,+)}=2.$

For $n\ge 3,$ we have
\begin{align}
  &\frac{\mathcal P_2^{n}(1,n,+)}{\mathcal P_2^{n+1}(1,n,+)}=\frac{P(Y'\text{ hits }[n,\infty)\text{ at }n\text{ before } [0,1]|Y_0'=2)}{P(Y'\text{ hits }[n,\infty)\text{ at }n+1\text{ before } [0,1]|Y_0'=2)}\no\\
  &\quad =\frac{P(Y'\text{ hits }[n,\infty)\text{ at }n\text{ before } [0,1], E_n|Y_0'=2)}{P(Y'\text{ hits }[n,\infty)\text{ at }n+1\text{ before } [0,1], E_n|Y_0'=2)}\no\\
  &\quad= \frac{\tilde{P}(T_n=n|Y_0'=2)}{\tilde{P}(T_n=n+1|Y_0'=2)}=\frac{h_n(1)}{h_n(2)}.\no
\end{align}
As a result, it follows from Step 4 that $\lim_{n\rto}\frac{\mathcal P_2^{n}(1,n,+)}{\mathcal P_2^{n+1}(1,n,+)}=\lim_{n\rto}\frac{h_n(1)}{h_{n}(2)}=2.$ The lemma is proved. \qed

\subsection{Connection between $\mathbf e_1N_2\cdots N_n\mathbf e_1$ and $\xi_2\cdots \xi_n$} \label{mt3}

We can now deduce from \eqref{xto}, \eqref{esb} and \eqref{rpk} that if $p_i=1/3\pm r_i,i\ge 2,$  then
\begin{align}\label{to}
   \mathcal P_2^n(1,n,+)\sim \frac{2}{3}\frac{1}{1+\sum_{i=2}^{n-1}\xi_{2}\cdots \xi_i}, \text{ as }n\rto.
\end{align}
With \eqref{to} in hands, to characterize the asymptotics of $P(M=n,D<\infty),$ we need to estimate further $\mathcal P_n(1,n+1,-),$ whose lower and upper bounds are given in \eqref{esa}. But the upper bound in \eqref{esa} is approximately twice as much as the lower bound, so that \eqref{esa} is not enough for us to get the accurate limit behavior of $\mathcal P_n(1,n+1,-).$ To deal with this difficulty, we come back to the product of nonnegative matrices.

By an argument similar to Lemma \ref{esy}, we have
\begin{align}
  \mathcal P_k(m,n,+)=\frac{\sum_{s=m+1}^{k}\mathbf e_1 N_s\cdots N_{n-1}\mathbf e_1^t}{1+\sum_{s=m+1}^{n-1}\mathbf e_1 N_s\cdots N_{n-1}\mathbf e_1^t},\no
\end{align}
from which, we get for $n\ge 2,$
\begin{align}
  &\mathcal P_n(1,n+1,-)=\frac{1}{1+\sum_{s=2}^{n}\mathbf e_1 N_s\cdots N_{n}\mathbf e_1^t}=\frac{1}{\sum_{s=2}^{n+1}\mathbf e_1 N_s\cdots N_{n}\mathbf e_1^t},\label{smd}\\
  &\mathcal P_2(1,n,+)=\frac{\mathbf e_1 N_2\cdots N_{n-1}\mathbf e_1^t}{1+\sum_{s=2}^{n-1}\mathbf e_1 N_s\cdots N_{n-1}\mathbf e_1^t}=\frac{\mathbf e_1 N_2\cdots N_{n-1}\mathbf e_1^t}{\sum_{s=2}^{n}\mathbf e_1 N_s\cdots N_{n-1}\mathbf e_1^t}.\label{sma}
\end{align}
 Although the rightmost-hand side of \eqref{smd} is written in terms of the product of nonnegative matrices $N_k,k\ge2,$  we can not use Theorem \ref{amij} to give an estimate because the summand in denominator involves $n.$ Thus we turn back to the continued fraction.
 \begin{lemma}\label{mc} If $p_i=1/3\pm r_i,i\ge2,$ then
 \begin{align}
 &(\xi_2\cdots\xi_n)^{-1}\sim c \mathbf e_1 N_2\cdots N_{n}\mathbf e_1^t\label{xen},\text{ as } n\rto,\\
   &\mathcal P_n(1,n+1,-)\sim c\frac{\xi_2\cdots\xi_n}{\sum_{s=2}^{n+1}\xi_2\cdots \xi_{s-1}}, \text{ as }n\rto.\label{pce}
 \end{align}
 \end{lemma}
 \proof For $2\le s\le n+1,$ set $$y_{s,n}:=\mathbf e_1 N_s\cdots N_{n}\mathbf e_1^t \text{ and }\xi_{s,n}:=\frac{y_{s+1,n}}{y_{s,n}}.$$ Noting that the empty product equals identity, thus $y_{n+1,n}=\mathbf e_1I\mathbf e_1=1.$ Therefore, we have
 \begin{align}\xi_{s,n}^{-1}\cdots \xi_{n,n}^{-1}=y_{s,n}=\mathbf e_1 N_s\cdots N_{n}\mathbf e_1^t.\label{cpn}\end{align}
Substituting \eqref{cpn} into \eqref{smd}, we obtain
\begin{align}\label{fpc}
  \mathcal P_n(1,n+1,-)=\frac{1}{\sum_{s=2}^{n+1}\xi_{s,n}^{-1}\cdots \xi_{n,n}^{-1}}=\frac{\xi_{2,n}\cdots \xi_{n,n}}{\sum_{s=2}^{n+1}\xi_{2,n}\cdots \xi_{s-1,n}}.
\end{align}
If we can show
\begin{align}
 & \xi_{2,n}\cdots \xi_{n,n}\sim c\xi_{2}\cdots \xi_{n},\label{xnp}\\
 & \sum_{s=2}^{n+1}\xi_{2,n}\cdots \xi_{s-1,n}\sim \sum_{s=2}^{n+1}\xi_{2}\cdots \xi_{s-1}, \label{sxn}
\end{align}
as $n\rto,$ then  \eqref{pce} is a consequence of \eqref{fpc}.

We prove first \eqref{xnp}. To this end, note that for $2\le s\le n,$
\begin{align}
  \xi_{s,n}&=\frac{y_{s+1,n}}{y_{s,n}}=\frac{\mathbf e_1 N_{s+1}\cdots N_{n}\mathbf e_1^t}{\mathbf e_1 N_s\cdots N_{n}\mathbf e_1^t}\no\\
  &=\frac{\mathbf e_1 N_{s+1}\cdots N_{n}\mathbf e_1^t}{\theta_s(\mathbf e_1+\mathbf e_2) N_{s+1}\cdots N_{n}\mathbf e_1^t}=\frac{\theta_s^{-1}}{1+\frac{\mathbf e_2 N_{s+1}\cdots N_{n}\mathbf e_1^t}{\mathbf e_1 N_{s+1}\cdots N_{n}\mathbf e_1^t}}.\no
  \end{align}
A backward induction shows that
$$\frac{\mathbf e_2 N_{s+1}\cdots N_{n}\mathbf e_1^t}{\mathbf e_1 N_{s+1}\cdots N_{n}\mathbf e_1^t}= \frac{\theta_{s+1}^{-1}}{1}\begin{array}{c}
                                \\
                               +
                             \end{array}\frac{\theta_{s+2}^{-1}}{1}\begin{array}{c}
                                \\
                               +\cdots
                             \end{array}\begin{array}{c}
                                \\
                               +
                             \end{array}\frac{\theta_n^{-1}}{1}.$$
Therefore we have
\begin{align}\xi_{s,n}=\frac{\theta_{s}^{-1}}{1}\begin{array}{c}
                                \\
                               +
                             \end{array}\frac{\theta_{s+1}^{-1}}{1}\begin{array}{c}
                                \\
                               +\cdots
                             \end{array}\begin{array}{c}
                                \\
                               +
                             \end{array}\frac{\theta_n^{-1}}{1}.\label{xs}\end{align}
                            Comparing \eqref{xs} with  \eqref{tail} and \eqref{dx}, by  theory of convergence  of limit periodic continued fractions, we have $\lim_{n\rto}\xi_{s,n}=\xi_s.$
By Lemma \ref{dfr}, the matrices $N_k^t,k\ge 2$ satisfies Condition (B2)$_c.$ Thus, applying Theorem \ref{amij}, we get
\begin{align}
\mathbf e_1 N_s\cdots N_{n}\mathbf e_1^t=\mathbf e_1 N_n^t\cdots N_{s}^t\mathbf e_1^t\sim c\varrho(N_2)\cdots \varrho(N_n).\no
  \end{align}
  Assume now $p_i=1/3\pm r_i,i\ge2.$ Then by \eqref{rs}, \eqref{sia} and \eqref{ra}, we have
    $$\varrho(N_2)\cdots \varrho(N_n)\sim  c\z(n\log n\cdots \log_{K-2}n(\log_{K-1} n)^B\y)^{\pm 1}, $$ and consequently
  $$\mathbf e_1 N_2\cdots N_{n}\mathbf e_1^t\sim   c\z(n\log n\cdots \log_{K-2}n(\log_{K-1} n)^B\y)^{\pm 1} \text{ as }n\rto.$$
  On the other hand, it follows from \eqref{xir} and \eqref{xa} that
  \begin{align}
    \xi_2\cdots\xi_n\sim c \z(n\log n\cdots \log_{K-2}n(\log_{K-1} n)^B\y)^{\mp 1} \text{ as }n\rto.\no
  \end{align}
  Therefore, we have
  \begin{equation*}
    (\xi_2\cdots\xi_n)^{-1}\sim c \mathbf e_1 N_2\cdots N_{n}\mathbf e_1^t, \text{ as }n\rto,
  \end{equation*}
  which proves \eqref{xen}.
  Taking \eqref{xen} and \eqref{cpn} together, we get \eqref{xnp}.

  Next, we proceed to prove \eqref{sxn}. For this purpose, noticing that by \eqref{sma} and \eqref{cpn},
    \begin{align}
\mathcal P_2(1,n+1,+)&=\frac{\mathbf e_1 N_2\cdots N_{n}\mathbf e_1^t}{\sum_{s=2}^{n+1}\mathbf e_1 N_s\cdots N_{n}\mathbf e_1^t}=\frac{\xi_{2,n}^{-1}\cdots \xi_{n,n}^{-1}}{\sum_{s=2}^{n+1}\xi_{s,n}^{-1}\cdots \xi_{n,n}^{-1}}\no\\
    &=\frac{1}{\sum_{s=2}^{n+1}\xi_{2,n}\cdots \xi_{s-1,n}}\no
\end{align}
 which together with \eqref{esb} yields that
 \begin{align}
      \frac{1}{1+\sum_{s=2}^{n+1}\xi_{2}\cdots \xi_s}\le \frac{1}{\sum_{s=2}^{n+1}\xi_{2,n}\cdots \xi_{s-1,n}}\le \frac{1}{1+\sum_{s=2}^{n}\xi_{2}\cdots \xi_s}.\no
 \end{align}
  As a consequence,
  \begin{align}
      1-\frac{\xi_{2}\cdots \xi_{n+1}}{1+\sum_{s=2}^{n+1}\xi_{2}\cdots \xi_s}\le \frac{\sum_{s=2}^{n+1}\xi_{2}\cdots \xi_{s-1}}{\sum_{s=2}^{n+1}\xi_{2,n}\cdots \xi_{s-1,n}}\le 1.\no
 \end{align}
 Therefore, using \eqref{xto} we get
 $$\lim_{n\rto}\frac{\sum_{s=2}^{n+1}\xi_{2}\cdots \xi_{s-1}}{\sum_{s=2}^{n+1}\xi_{2,n}\cdots \xi_{s-1,n}}=1.$$
 Then \eqref{sxn} is proved and so is the lemma. \qed

\subsection{Proof of Theorem \ref{myp}}\label{mt4}
 We now give the proof of Theorem \ref{myp}.
 Suppose that $p_i=1/3\pm r_i,i\ge 2.$ Then by \eqref{to} and \eqref{pce}, we have
 \begin{align}
   P(M=n&, D<\infty)=\mathcal P_2^n(1,n,+)\mathcal P_n(1,n+1,-)\no\\
&\sim c\frac{1}{\sum_{s=2}^{n}\xi_{2}\cdots \xi_{s-1}}\times \frac{\xi_2\cdots\xi_n}{\sum_{s=2}^{n+1}\xi_2\cdots \xi_{s-1}}, \text{ as }n\rto.\no
 \end{align}
 Furthermore, it follows from Lemma \ref{xie} that
 $$\xi_2\cdots \xi_n\sim  c\z(n\log n\cdots \log_{K-2}n(\log_{K-1} n)^B\y)^{\mp 1},\text{ as }n\rto.$$
With the above facts in hands, once again, the proof of Theorem \ref{myp} is just a step-by-step repetition   of that of \cite[Theorem 1]{wh}. \qed

\section*{Appendix:  Recurrence criteria }
We present here the criteria of recurrence, positive recurrence and transience for both the chain $Y$ and its adjoint chain $Y',$ since they are helpful to understand the various decay speeds in Theorem \ref{mxt}  and Theorem \ref{myp}.

Let $\theta_k,N_k,k\ge 2$ be those defined in \eqref{tn} and $\xi_k,k\ge2$ be the one in \eqref{dx}. Then
 as a direct consequence of Lemma \ref{esy} and Lemma \ref{ecpb},  we  get the transience criteria of both  $Y$ and $Y'.$
\begin{corollary}\label{tc} {\rm (i)} The (2,1) random walk $Y$ is transient or recurrent according as $\sum_{s=2}^{\infty}\mathbf e_1 N_s\cdots N_{2}\mathbf e_1^t <\infty $ $\text{ or }=\infty.$ {\rm(ii)} The (1,2) random walk $Y'$ is transient or recurrent according as $\sum_{s=2}^\infty \xi_2\cdots \xi_s<\infty\text{ or }=\infty.$
\end{corollary}
The following recurrence criteria  of $Y$ and its adjoint chain $Y'$ can be found in  \cite[page 204]{dey98}.
\begin{lemma}\label{yyp}
  {\rm (i)} The chain $Y$ is positive recurrent if and only if its adjoint chain $Y'$ is transient and vice versa.
     {\rm (ii)} Both adjoint chains $Y$ and $Y'$ are null recurrent simultaneously.
\end{lemma}
 Let $r_i$ be the one  in (\ref{ri}). If $p_i=1/3\pm r_i,i\ge2,$ then by Lemma \ref{xie} and Lemma \ref{mc}, we have
 \begin{align}
 (\xi_2&\cdots\xi_n)^{-1}\sim c \mathbf e_1 N_2\cdots N_{n}\mathbf e_1^t\no\\
& \sim c\z(n\log n\cdots \log_{K-2}n(\log_{K-1} n)^B\y)^{\pm 1}\text{ as } n\rto.\label{xne}
 \end{align}
 Combining Corollary \ref{tc} and Lemma \ref{yyp} with \eqref{xne}, we have the following result.
 \begin{corollary}\label{crpt}
{\rm (i)} For $K=1,$ if $q_i=\frac{2}{3}+r_i,\ i\ge 2$(or $p_i=\frac{1}{3}-r_i,\ i\ge 2$),
then
\begin{align}
  B>1&\Rightarrow  Y \text{ is transient  and }Y' \text{ is positive recurrent;}\no\\
  B<-1&\Rightarrow  Y' \text{ is transient  and }Y \text{ is positive recurrent;}\no\\
B\in [-1,1]&\Rightarrow \text{both }Y\text{ and }Y' \text{ are null recurrent.}\no
\end{align}
 {\rm (ii)}  For $K\ge 2,$ if $q_i=\frac{2}{3}+r_i,\ i\ge 2,$ then
   \begin{align}
  &B>1\Rightarrow Y \text{ is transient  and }Y' \text{ is positive recurrent;}\no\\
  &B\le1 \Rightarrow \text{both }Y\text{ and }Y' \text{ are null recurrent.}\no
\end{align}
{\rm (iii)} For $K\ge 2,$ if $q_i=\frac{2}{3}-r_i,\ i\ge 2,$ then
   \begin{align}
  &B>1\Rightarrow Y' \text{ is transient  and }Y \text{ is positive recurrent;}\no\\
  &B\le1 \Rightarrow \text{both }Y\text{ and }Y' \text{ are null recurrent.}\no
\end{align}
\end{corollary}
\proof
Here, we give only the proof of part (iii), since  the others can be proved similarly.
Fix $K= 2,3,...$ and let $q_i=\frac{2}{3}-r_i,\ i\ge 2.$ Then, by \eqref{xne},  \begin{equation}\xi_2\cdots\xi_n\sim \frac{1}{\mathbf e_1 N_n\cdots N_{2}\mathbf e_1^t}\sim c\frac{1}{n\log n\cdots \log_{K-2}n(\log_{K-1} n)^B}\label{cd}\end{equation}
 as  $n\rto.$
If $B>1,$ then $\sum_{s=2}^\infty\xi_{2}\cdots\xi_s<\infty.$ Consequently, it follows from Corollary \ref{tc} that $Y'$ is transient. Therefore, by Lemma \ref{yyp}, $Y$ must be positive recurrent.

If $B\le1,$ then $\sum_{s=2}^\infty\xi_{2}\cdots\xi_s=\infty.$  Thus, by Corollary \ref{tc}, $Y'$ is recurrent. We claim that $Y'$ must be null recurrent. Indeed, if we suppose  conversely $Y'$ is positive recurrent, then by Lemma \ref{yyp}, its adjoint chain $Y$ must be transient. So it follows from Corollary \ref{tc} that $\sum_{s=2}^{\infty}\mathbf e_1 N_s\cdots N_{2}\mathbf e_1^t <\infty$ which contradicts  \eqref{cd}. Therefore we can conclude that $Y'$ is null recurrent and so is its adjoint chain  $Y.$ The proof of part (iii) is finished. \qed

\vspace{0.1cm}

\noindent{\large{\bf \Large Acknowledgements:}} The authors would like to thank Prof. W. M. Hong for introducing to us the Lamperti problem and Prof. Y. Chow for some useful comments. This project is supported by National
Natural Science Foundation of China (Grant No. 11501008; 11601494).

\end{document}